\documentclass[a4paper,12pt]{amsart}

\usepackage[none]{hyphenat}
\usepackage{a4wide,hyperref,amsmath,amssymb,enumitem, stmaryrd}
\usepackage[normalem]{ulem}
\usepackage{tikz}
\usetikzlibrary{decorations.pathreplacing,calligraphy}
\usetikzlibrary{shapes}
\tikzset{
			inner sep=1pt,semithick,
			vertex/.style={circle,draw,fill,minimum size= 1pt},
			vertexb/.style={rectangle,draw,fill, minimum size = 3pt},
			vertexd/.style={rectangle,draw,fill=white, minimum size = 3pt},			
			vertexc/.style={circle,draw,fill=white,minimum size= 1pt},			
			thickedge/.style={line width=1pt},
			rededge/.style={line width=2pt, black},
			blueedge/.style={line width=1pt, black},
			wideyellow/.style={line width=3pt,yellow},						
			font=\tiny
}

\usepackage{lineno}

\newcommand{\X}{\mathcal{X}}

\newcommand{\Y}{\mathcal{Y}}

\newif\ifdetails
\detailstrue
\newcommand{\DETAIL}[1]%
{\ifdetails\par\fbox{\begin{minipage}{0.9\linewidth}\textit{Detail:}
      #1\end{minipage}}\par\fi}
\newcommand{\TODO}[1]%
{\ifdetails\par\fbox{\begin{minipage}{0.9\linewidth}\textbf{TODO:}
      #1\end{minipage}}\par\fi}

\newtheorem{lemma}{Lemma}

\newtheorem{theorem}[lemma]{Theorem}

\newtheorem{definition}{Definition}

\newtheorem{problem}{Problem}
\newtheorem{claim}{Claim}

\newcommand{\lcr}{\overline{\operatorname{cr}}}
\newcommand{\Cr}{\operatorname{cr}}

\newcommand{\old}[1]{{}}

\makeatletter
\DeclareRobustCommand{\cev}[1]{%
  {\mathpalette\do@cev{#1}}%
}
\newcommand{\do@cev}[2]{%
  \vbox{\offinterlineskip
    \sbox\z@{$\m@th#1 x$}%
    \ialign{##\cr
      \hidewidth\reflectbox{$\m@th#1\vec{}\mkern4mu$}\hidewidth\cr
      \noalign{\kern-\ht\z@}
      $\m@th#1#2$\cr
    }%
  }%
}
\makeatother

\title[A Quantitative Kuratowski Theorem]{A Quantitative Kuratowski Theorem}

\author[\'E. Czabarka, A. Helm, G. T\'oth]{\'Eva Czabarka, Alec Helm, G\'eza T\'oth}
\address{\'Eva Czabarka, Alec Helm\\ Department of Mathematics \\ University of South Carolina \\ Columbia, SC 29208 \\ USA\\
G\'eza T\'oth\\Alfr\'ed R\'enyi Institute of Mathematics
\\Budapest, Hungary}

\thanks{G\'eza T\'oth was supported by  ERC 
Advanced Grant `GeoScape' No.\ 882971 and by the National Research, Development and Innovation Office, NKFIH, ADVANCED 152590.}

\date{\today}

\begin{document}


\begin{abstract}
Suppose that a graph $G$ contains exactly $k$ different subdivisions of a $K_5$ or $K_{3,3}$.
We prove that in this case the crossing number of $G$ is at most $\binom{3k+7}{2}$ (and, when $k\le 2$, it is at most $k$). 
\end{abstract}

\maketitle


\section{Introduction}

In this paper a graph is finite, may have multiple edges but no loops. 
As usual, $|G|$ is the number of vertices of $G$.
In a {\em drawing} 
$\mathcal{D}$
of a graph $G$, vertices are represented by distinct points in the plane, edges are represented by curves, connecting the corresponding points. We call these points and curves also vertices and edges, respectively.
We assume that an edge does not contain a vertex in its interior, two edges have finitely many common points and each of them is a common endpoint or a proper crossing. Finally, we assume that no three edges cross at the same point.
The {\it crossing number} $\Cr(G,{\mathcal D})$ of a drawing $\mathcal D$ of a graph $G$
is the number of edge crossings, 
and the {\it crossing number} $\Cr(G)$ of a graph $G$
is the minimum crossing number over all of its drawings. 
A {\it subdivision} of a graph is obtained by replacing some of its edges with internally vertex disjoint paths. 
Two graphs are topologically isomorphic if they are subdivisions of the same graph.

It is easy to see that the crossing number of a graph is the same as the crossing number of any of its subdivisions. A graph $G$ is planar if and only if $\Cr(G)=0$. By  
the celebrated Kuratowski Theorem, \cite{kuratowski} $\Cr(G)>0$,  that is to say that $G$ is non-planar, if and only if it contains a subdivision of either $K_5$ (the complete graph on $5$ vertices) or 
$K_{3,3}$ (the complete bipartite graph with $3$ vertices in each partition class) as a subgraph. 
We will refer to subdivisions of $K_5$ and $K_{3,3}$ 
as {\it cross-inducing graphs}. 

A graph $G$ is {\it $m$-crossing-critical} if $\Cr(G)\geq m$ but the crossing number of any proper subgraph of $G$ is below $m$. It follows from 
the Kuratowski Theorem that
$1$-crossing-critical graphs are exactly the subdivisions of $K_5$ and $K_{3,3}$. Kochol \cite{kochol} gave an infinite family of $3$-connected $k$-crossing-critical graphs for any $k\ge 2$.
Apart from finitely many possible exceptions, the $2$-crossing-critical graphs were characterized by Bokal, Oporowski, Richter and Salazar \cite{2cc}. 
Characterizing $k$-crossing-critical graphs for $k>2$ seems hopeless at the moment.
See \cite{bdc}, \cite{bt22}, and \cite{curious} for some recent results on crossing-critical graphs.

Let $k$ be an arbitrary positive integer. It is possible that a graph $G$ 
contains $k$ different cross-inducing subgraphs, and still $\Cr(G)=1$. 
As an example, take a $K_5$ and replace one of its edges by $k$ internally disjoint paths. 
But what is the {\em maximum} crossing number of a graph which contains exactly $k$ different (not necessarily disjoint) 
cross-inducing subgraphs?
Let $f(k)$ denote this quantity, if this maximum exists, and $f(k)=\infty$ otherwise. 

Obviously, $f$ is an increasing function and $f(k)\ge k$ as shown by a disjoint union of $k$ $K_5$'s. 
This is the best lower bound we have for $f$. 
The function 
$f$ is super-additive as shown by taking disjoint unions of graphs 
$G_1$ and $G_2$ that have $t_1$ and $t_2$ cross-inducing subgraphs and have crossing number $f(t_1)$ and $f(t_2)$, respectively. 
One might suspect that at least $k$ cross-inducing subgraphs are needed for a graph to have crossing number $k$, meaning that $f(k) = k$, but 
it is not immediately obvious that $f(k)$ is even finite for every $k$.

The Kuratowski Theorem implies that $f(0)=0$. 
In this paper we strengthen this result to the following.
\begin{theorem}\label{th:geza} 
Let $G$ be a graph with $k$ cross-inducing subgraphs.
If $k\le 2$, $G$ has a drawing with at most $k$ crossings.
In particular, for $k\leq 2$ we have $f(k)=k$.
\end{theorem}

Note that $f(k)\le m$ is equivalent with the statement that every $(m+1)$-crossing-critical graph has at least $k+1$ cross-inducing subgraphs,
therefore this gives that the $2$-crossing critical graphs almost characterized in \cite{2cc} each have at least $2$ cross-inducing subgraphs.

Our main result states that $f(k)$ is indeed finite for every $k$:
\begin{theorem}\label{th:existence} If $G$ is a graph with at most $k$ cross-inducing subgraphs, then it has a drawing with at most
$\binom{3k+7}{2} =\frac{9}{2}k^2+\frac{39}{2}k+21$ crossings.
Consequently, 
$k\le f(k) \le \binom{3k+7}{2}$.
\end{theorem}

\begin{problem} Is it true that $f(k)=k$ for every $k\ge 1$? 
\end{problem}

The proof of Theorem \ref{th:existence} is quite long, it is divided into several lemmas. 
Surprisingly, we do not have a shorter proof for any weaker bound.



A simple graph is a graph with no multiple edges. The F\'ary--Wagner theorem \cite{fary}, \cite{wagner}
states that any simple, planar graph can be drawn without crossings with straight line segments as edges, 
i.e. the 
cross-inducing subgraphs with respect to straight line drawings are still the subdivisions of $K_{3,3}$ and $K_5$. 
The rectilinear crossing number $\lcr(G)$ of a simple graph $G$ is the minimum crossing number over its straight line drawings. 
Clearly, for any simple graph $G$ we have $\lcr(G)\ge\Cr(G)$. 
Bienstock and Dean \cite{bienstock} proved a generalization of the F\'ary--Wagner theorem. They proved that 
if $\Cr(G)\le 3$ for a graph $G$, then $\Cr(G)=\lcr(G)$. On the other hand, 
for any $m \ge  4$, they constructed
a class of graphs with $\Cr(G)=m$ and $\lcr(G)$ arbitrarily large.

Let $\bar{f}(k)$ denote the {\em maximum} rectilinear crossing number of a simple 
graph which contains exactly $k$ different 
cross-inducing subgraphs. (If the maximum exists, $\bar{f}(k)=\infty$ otherwise.)
As Bienstock and Dean \cite{bienstock} showed that the crossing number and the rectilinear crossing number are the same when the crossing number is at most $3$, Theorem~\ref{th:geza} yields that
$\bar{f}(k)=f(k)=k$ for $k\le 2$.

We have 
$\lcr(G)\ge\Cr(G)$ for any simple graph, but it does not imply automatically that 
$\bar{f}(k)\ge {f}(k)$ since for $f(k)$ we also consider multigraphs. However, 
any multigraph $G$ can be converted to a simple graph $G^*$
by subdividing all but one parallel edges between any two vertices, and $\lcr(G^*)\ge \Cr(G^*)=\Cr(G)$,
 $\bar{f}(k)\ge {f}(k)$. 
This implies that  $\bar{f}(k)\ge k$, and it is the best lower bound we have for 
 $\bar{f}(k)$.

\begin{problem} Is $\bar{f}(k)$ finite for every $k\ge 3$?
\end{problem}

The organization of the remainder of this paper is the following: 
Section~\ref{sec:prelim} gives the necessary notations and definitions, and provides some equivalent statements with our theorems. 
Section~\ref{sec:lemmas} gives basic results
on how adding new edges to a $K_5$ or $K_{3,3}$
adds new cross-inducing subgraphs, while
Section~\ref{sec:structure}  is concentrated on 2-connected graphs, and reduces the proof of our theorems to finding specific types of drawings in certain 2-connected graphs.
Section~\ref{sec:k33Linear} gives some technical edge-counting necessary for our bound in Theorem~\ref{th:existence}, and we give the proofs of Theorem~\ref{th:geza} and~\ref{th:existence} in Section~\ref{sec:mainproof}.

\section{Preliminaries}\label{sec:prelim}

As is customary, we use the notation
$[n]=\{1,2,\ldots,n\}$. 

\begin{definition}
 Let $G$ be a graph. A vertex $v$ is a \emph{subdividing vertex}, 
if it has degree $2$, otherwise it is a \emph{main vertex}.

Let $H$ be a subgraph of the graph $G$. We say a vertex $v$ of $G$ is $\emph{$H$-main}$ if it is a main vertex of $H$, and it is  $\emph{$H$-subdividing}$ if it is a subdividing vertex of $H$. 
We say that a pair of vertices $u,v$ are \emph{$H$-cooperative} if $H$ is a subdivided $K_{3,3}$, and $u$ and $v$ are $H$-main vertices in the same main vertex class.
\end{definition}

\smallskip

\begin{definition}
Given a graph $G$, an \emph{edge-path} of $G$ is a path $P$ of length at least $1$ such that the endvertices of $P$ are main vertices, and the internal vertices of $P$ (if any) are subdividing vertices.

If $H$ is a subgraph of $G$, an \emph{$H$-edge path} is an edge-path of $H$.

 If $H$, $G$ are graphs ($H$ is not necessarily a subgraph of $G$), then a path $P$ is called a 
\emph{$(G,H)$-path} if $P$ is a path in $G$ and no edge or 
internal vertex of $P$ is in $H$. 
\end{definition} 

\smallskip

\begin{definition}
Let $G$ be a graph with cross-inducing subgraphs $X_1,...,X_k$. We set $\X_G = \bigcup_{i\in[k]}X_i$, ie the union of all cross-inducing subgraphs of $G$.  $\Y_G$ is the edge-induced subgraph with edge set $E(G)-E(\X_G)$. 
\end{definition}

\smallskip

Notice that for any graph $G$, $\Cr(\Y_G) = 0$, since $\Y_G$ has no cross-inducing subgraphs. 
Clearly, $\Cr(\X_G)\le\Cr(G)$, and
\u{S}ir\'a\u{n} \cite{siran} gave an example of a graph $G$ such that $\Cr(\X_G)<\Cr(G)$, which \cite{curious} strengthened to an example of a graph $G$ where in every optimal drawing of $G$ some edge of $\Y_G$ is crossed.

Our proof will use the following strategy: 
given a graph $G$ with $k$ cross-inducing subgraphs, in order to provide 
an upper bound $g(k)$ for $f(k)$, we create a drawing of $G$
with at most $g(k)$ crossings, in which only edges of $\X_G$ are crossed. 
The example in \cite{curious} shows that in general it is possible that in every 
optimal drawing of $G$, some edges of $\Y_G$ are crossed. However, in our case, 
we can assume without loss of generality, that edges of $\Y_G$ do not participate in any crossing.
\begin{lemma}\label{lm:noYcross}
Let $k,M\in\mathbb{N}$. The following are equivalent
\begin{enumerate}[label={\upshape (\arabic*)}]
\item\label{case:generaldraw} Every graph $G$ with $k$ cross-inducing subgraphs has a drawing with at most $M$ crossings.
\item\label{case:reduceddraw} Every graph $G$ with  $k$ cross-inducing subgraphs has a drawing with at most $M$ crossings in which edges of $\Y_G$ are not crossed.
\end{enumerate}
\end{lemma}
\begin{proof}
\ref{case:reduceddraw} implies~\ref{case:generaldraw}, so assume ~\ref{case:generaldraw} holds for $k,M$.
Assume $G$ is graph with at most $k$ cross-inducing subgraphs. Let $G^*$ be the graph obtained from $G$ by replacing every edge of $\Y_G$ by $M+1$ parallel edges. Then $G^*$ and $G$ have the same cross-inducing subgraphs, and $\X_G=\X_{G^*}$. 
Select an optimal drawing $\mathcal{D}$ of $G^*$. By~\ref{case:generaldraw},
$\Cr(\mathcal{D},G^*)\le M$. In any optimal drawing, parallel edges are crossed the same number of times. Since $M+1>M$, no edge
of $\Y_{G^*}$ is crossed. Thus, $\mathcal{D}$ contains a drawing of $G$ with at most $M$ crossings in which no edge of $\Y_G$ cross, i.e. ~\ref{case:reduceddraw} holds.
\end{proof}

Let $S^2$ be the surface of the sphere and $p$ be a point (the ``northern pole") of the sphere, then $s^2\setminus\{p\}$ is mapped to the plane $R^2$ by a stereographic projection. 
A drawing of a graph in the plane $\mathbb{R}^2$ can be mapped to a drawing of the graph on the surface $S^2\setminus\{p\}$ of the sphere through an inverse of a stereographic projection, and using a different stereographic projection it can be mapped back into another drawing.

\begin{definition}
The drawings $\mathcal{D}$, $\mathcal{D}'$ 
of a graph $G$ are equivalent, if
there is a homeomorphism $\phi$ of a sphere and stereographic projections 
$\psi_1,\psi_2$, $\psi_i:S^2\setminus\{p_i\}\rightarrow\mathbb{R}^2$ such $\psi_2 \circ \phi \circ \psi_1^{-1}$
maps $\mathcal{D}$ to $\mathcal{D}'$.
\end{definition}

We will use the following easy and well-known statements.  
\begin{lemma}\label{lm:obvious} Let $G_1,G_2$ be two graphs with drawings $\mathcal{D}_1,\mathcal{D}_2$, respectively. 
If $G_1\cap G_2$ is empty or 
a single vertex, then $G_1\cup G_2$ has a drawing that induces a drawing of each $G_i$ 
equivalent with $\mathcal{D}_i$ in which no edge of $G_1$ crosses an edge of $G_2$; 
in particular $\Cr(G_1\cup G_2)=\Cr(G_1)+\Cr(G_2)$.
If $G_1\cap G_2$ is two vertices connected by an edge $e$, and for each $i\in\{1,2\}$ 
the drawing $\mathcal{D}_i$ has no crossings involving $e$, 
then $G_1\cup G_2$ has a drawing that induces a drawing of each $G_i$ equivalent with 
$\mathcal{D}_i$ respectively, such that no edge of $G_1$ crosses an edge of $G_2$; 
in particular $\Cr(G_1\cup G_2)\le \Cr( \mathcal(G_1,{D}_1)+\Cr( G_2, \mathcal{D}_2))$.\end{lemma}

The next lemma  shows that we only need to prove Theorems~\ref{th:geza} and~\ref{th:existence} 
for $2$-connected graphs and for drawings where only edges of $\X_G$ cross.

\begin{lemma}\label{lm:superadditive}
Let $g:\mathbb{N}\rightarrow\mathbb{N}$ be a superadditive function and $k\in\mathbb{N}$. The following are equivalent
\begin{enumerate}[label={\upshape (\arabic*)}]
    \item\label{case:reduce2plus} For all $k'\le k$ 
every $2$-connected graph $H$ with exactly $k'$ cross-inducing subgraphs 
has a drawing with at most $g(k')$ crossings in which only edges of $\X_H$ cross.
    \item \label{case:reduce1minus} For all $k'\le k$ every graph $G$ with exactly $k'$ cross-inducing subgraphs has a drawing with at most $g(k')$ crossings.
    \end{enumerate}
\end{lemma}

\begin{proof}
\ref{case:reduce1minus} together with Lemma~\ref{lm:noYcross} implies~\ref{case:reduce2plus}.
Assume~\ref{case:reduce2plus} holds.
Let $G$ be a graph with exactly $k'$ cross-inducing subgraphs for some $k'\le k$.
Let $B_1,\ldots,B_t$ be
the 
blocks of the graph $G$, and for $i\in[t]$ let $k_i$ be the number of cross-inducing subgraphs 
of $B_i$. In particular, $0\le k_i\le k'\le k$. 
Since $B_i\subseteq G$, $\bigcup_i\X_{B_i}\subseteq\X_G$, so $\sum_ik_i\le k'$. 
For every $i$, $B_i$ is either $2$-connected or $B_i$ is a $K_2$ -- in the latter case $\Cr(B_i)=0$.
If $B_i$ is $2$-connected, by~\ref{case:reduce2plus}
$B_i$ has a drawing $\mathcal{D}_i$ with at most $g(k_i)$ crossings in which only edges of $\X_{B_i}$ cross.
Lemma~\ref{lm:obvious} gives
a drawing of $G$ in which two edges cross precisely when they cross in some drawing
$\mathcal{D}_i$. 
Thus, we have a drawing of $G$ with at most 
$\sum_i g(k_i)\le g\left(\sum_i k_i\right)\le g(k')$ crossings
in which only edges of $\bigcup_i\X_{B_i}$
cross. (Note that we used that superadditive nonnegative functions must be increasing). As $\bigcup_i\X_{B_i}\subseteq\X_G$, the statement follows. 
\end{proof}


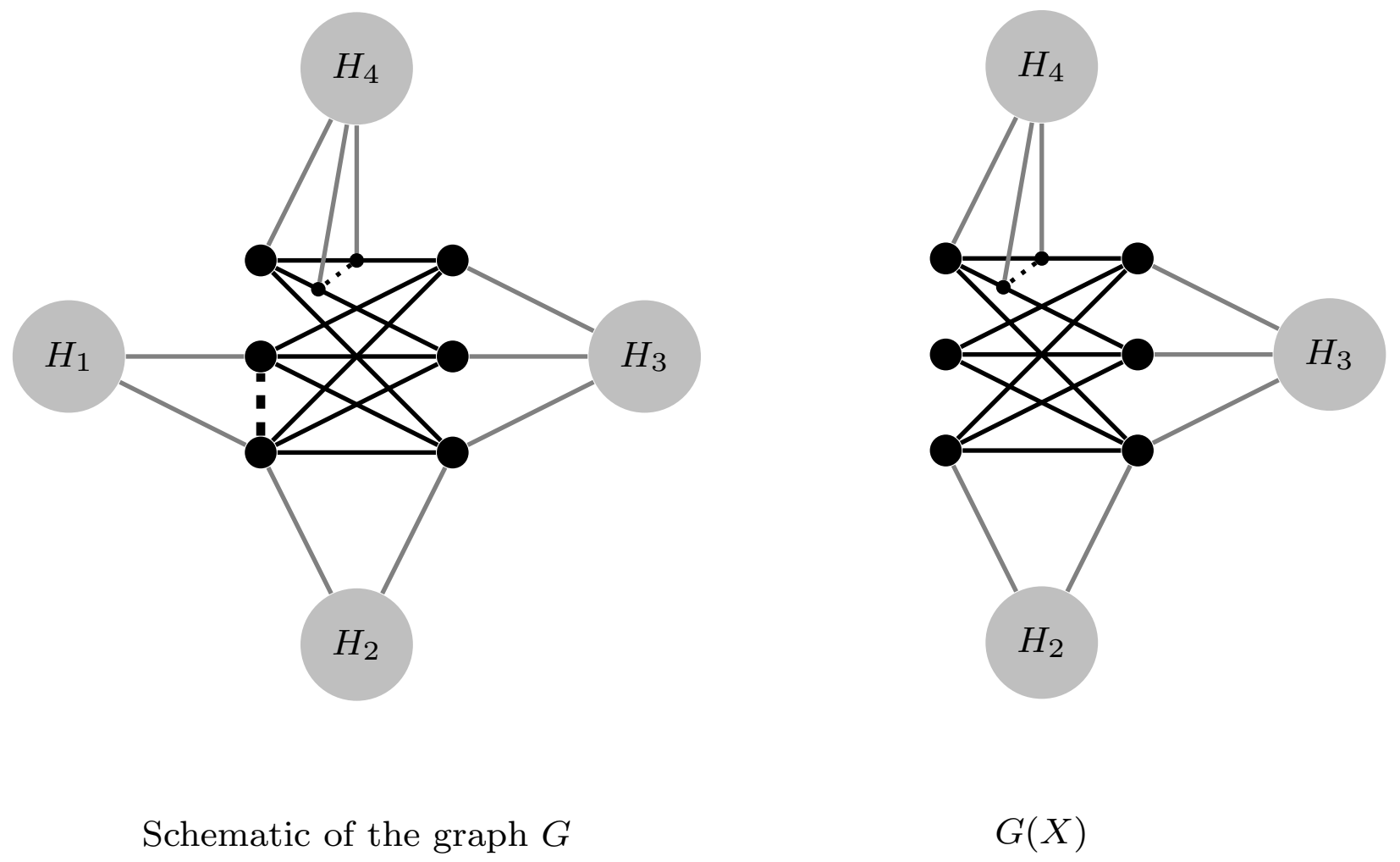
\begin{figure}[http]
\center
\begin{tikzpicture}[scale = 1.5]
\node[circle, fill=black, inner sep=2.5pt] (a1) at (1,0) {};
\node[vertex,black] (c1) at (.5,1) {};
\node[vertex,black] (c2) at (.3, .85) {};
\node[circle, fill=black, inner sep=2.5pt] (a2) at (1,.5) {};
\node[circle, fill=black, inner sep=2.5pt] (a3) at (1,1) {};
\node[circle, fill=black, inner sep=2.5pt] (b1) at (0,0) {};
\node[circle, fill=black, inner sep=2.5pt] (b2) at (0,.5) {};
\node[circle, fill=black, inner sep=2.5pt] (b3) at (0,1) {};
\draw[rededge, dashed] (b1)--(b2);
\draw[thickedge,dotted] (c1)--(c2);
\draw[thickedge,black] (a1)--(b1)--(a2)--(b2)--(a3)--(b3)--(a1);
\draw[thickedge,black] (a1)--(b2);
\draw[thickedge,black] (a2)--(b3);
\draw[thickedge,black] (a3)--(b1);
\node[circle,fill=gray!50, minimum size= 25pt] (h1) at (-1,.5) {$H_1$}; 
\node[circle,fill=gray!50, minimum size= 25pt] (h2) at (.5,-1) {$H_2$}; 
\node[circle,fill=gray!50, minimum size= 25pt] (h3) at (2,.5) {$H_3$}; 
\node[circle,fill=gray!50, minimum size= 25pt] (h4) at (.5,2) {$H_4$}; 
\draw[thickedge,gray] (b1)--(h1);
\draw[thickedge,gray](b2)--(h1);
\draw[thickedge,gray](h2)--(b1);
\draw[thickedge,gray](h2)--(a1);
\draw[thickedge,gray](h3)--(a1);
\draw[thickedge,gray](h3)--(a2);
\draw[thickedge, gray](h3)--(a3);
\draw[thickedge,gray](h4)--(c1);
\draw[thickedge,gray](h4)--(c2);
\draw[thickedge,gray](h4)--(b3);
\node at  (.5,-2) {Schematic of the graph $G$};
\end{tikzpicture}
\qquad\qquad
\begin{tikzpicture}[scale=1.5]

\node[circle, fill=black, inner sep=2.5pt] (a1) at (1,0) {};
\node[vertex,black] (c1) at (.5,1) {};
\node[vertex,black] (c2) at (.3, .85) {};
\node[circle, fill=black, inner sep=2.5pt] (a2) at (1,.5) {};
\node[circle, fill=black, inner sep=2.5pt] (a3) at (1,1) {};
\node[circle, fill=black, inner sep=2.5pt] (b1) at (0,0) {};
\node[circle, fill=black, inner sep=2.5pt] (b2) at (0,.5) {};
\node[circle, fill=black, inner sep=2.5pt] (b3) at (0,1) {};

\draw[thickedge,dotted] (c1)--(c2);
\draw[thickedge,black] (a1)--(b1)--(a2)--(b2)--(a3)--(b3)--(a1);
\draw[thickedge,black] (a1)--(b2);
\draw[thickedge,black] (a2)--(b3);
\draw[thickedge,black] (a3)--(b1);
\node[circle,fill=gray!50, minimum size= 25pt] (h2) at (.5,-1) {$H_2$}; 
\node[circle,fill=gray!50, minimum size= 25pt] (h3) at (2,.5) {$H_3$}; 
\node[circle,fill=gray!50, minimum size= 25pt] (h4) at (.5,2) {$H_4$}; 
\draw[thickedge,gray](h2)--(b1);
\draw[thickedge,gray](h2)--(a1);
\draw[thickedge,gray](h3)--(a1);
\draw[thickedge,gray](h3)--(a2);
\draw[thickedge, gray](h3)--(a3);
\draw[thickedge,gray](h4)--(c1);
\draw[thickedge,gray](h4)--(c2);
\draw[thickedge,gray](h4)--(b3);
\node at  (.5,-2) {$G(X)$};
\end{tikzpicture}
\caption{\small The schematic of the graph $G$ contains a $K_{3,3}$ subdivision $X$ (black vertices and black lines as edges, $X$-subdividing vertices are smaller), four components $H_1,\ldots,H_4$ of the graph $G-X$ (gray circles) connected to the graph by the gray edges, a dotted edge that is an element of $E_X^*$ and a thick dashed edge that is an element of $F_X^*$. $\mathbb{P}_1^X=\{H_1\}$, $\mathbb{P}_2^X=\{H_2,H_3,H_4\}$. }
\label{fig:G(X)def}
\end{figure}

 The proof will proceed by creating these drawings for $2$-connected graphs. 
 When all cross-inducing subgraphs are subdivided $K_{3,3}$-s, 
 we rely on
the following technical definition/notations; in fact, much of the work will 
involve determining properties of this structure. For illustration, see Figure~\ref{fig:G(X)def}.
\begin{definition}\label{def:k33structure}
Let $G$ be a graph and $X$  be a subdivided $K_{3,3}$. Then:
\begin{enumerate}[label={\upshape (\arabic*)}]
\item For  $X$-cooperative vertices $s,t$ , $\gamma_{st}$ denotes an edge not in $G$ that connects $s$ and $t$.
\item $E^*_X$ is the set of edges of $G$ connecting two $X$-cooperative vertices. 
\item $F^*_X$ is the set of edges of $G$ connecting two 
vertices of $X$ that are not $X$-cooperative  
(either at least one is a subdividing vertex, or they are in different main vertex classes).
\item $\mathbb{P}_X$ denotes the set of components of the graph $G-X$.
\item For a component $H\in\mathbb{P}_X$, $F_H$ is the set of edges 
that connect $H$ to $X$, $N_H$ is the set of vertices of 
$X$ that are endpoints of some edge in $F_H$ and $\widehat{H}=H\cup N_H\cup F_H$.
\item 
$\mathcal{P}_1^X=\left\{H\in\mathbb{P}_X: N_H=\{s_H,t_H\}\text{ where $s_H,t_H$ are $X$-cooperative}\right\}$,
$ 
\mathcal{P}_2^X=\mathbb{P}_X-\mathcal{P}_1^X
$.

\item $G(X)=X\cup F^*_X\cup\left(\bigcup_{H\in\mathcal{P}_2^X}\widehat{H}\right)$ and $m(X)=\sum_{H\in\mathcal{P}_2^X}|H|$.
\end{enumerate} 
\end{definition}

\section{Setting the stage}\label{sec:lemmas}

\begin{lemma}~\label{lm:onecross} Let $G$ be a graph and let $X$ be a cross-inducing subgraph of $G$. Assume there are two vertices $p\ne q$ of $X$ such that there is a $(G,X)$-path $P$ between $p$ and $q$. Then one of the following must hold:
\begin{enumerate}[label={\upshape (\roman*)}]
\item\label{case:one} $X$ is a subdivided $K_{3,3}$ and $p,q$ are $X$-cooperative
\item\label{case:two} $p$ and $q$ are (main or subdividing) vertices on the same $X$-edge-path $Q$. Replacing the $p$-$q$ sub-path  of $Q$ with $P$ we get another cross-inducing subgraph $Y$ topologically isomorphic to $X$ and the main vertices of $X$ and $Y$ are the same.
\item\label{case:twoprime} $q$ is $X$-subdividing, $p$ is $X$-main and not on the same edge-path as $q$. Then there is a cross-inducing subgraph $Y\subset X\cup P$ which is a subdivided $K_{3,3}$. It contains $P$ and its main vertices are in $X$.  
\item\label{case:three} Vertices $p$ and $q$ are both $X$-subdividing and are not on the same $X$-edge-path. Then $G$ has at least $3$ cross-inducing subgraphs in $X\cup P$.
\end{enumerate}
\end{lemma}

\begin{proof}
The conditions of the cases on $p$ and $q$ clearly cover all possibilities.
The statement of ~\ref{case:two} is obvious, and~\ref{case:one} does not require proof.
So assume that neither~\ref{case:one} nor the conditions in ~\ref{case:two} hold. Then at least one of $p,q$ is $X$-subdividing, and
$p,q$ are on different $X$-edge-paths. Due to the symmetries of $K_5$ and $K_{3,3}$
it is enough to check the cases on Figure~\ref{fig:geza1} when $X$ is a subdivided $K_5$, and on Figure~\ref{fig:geza2} when $X$ is a subdivided $K_{3,3}$. The statement follows.

 \begin{figure}[http]
\begin{tikzpicture} 
\node[vertexc] (a2) at (0,0) {$a_2$};
\node[vertexc] (a3) at (1,0) {$a_3$};
\node[vertexc] (a4) at (2,1) {$a_4$};
\node[vertexc] (a1) at (0,2) {$a_1$};
\node[vertexc] (a5) at (1,2) {$a_5$};
\node[vertexc] (p) at (.35,1.25) {$p$};
\node[vertexc] (q) at (0,.6) {$q$};
\draw[blueedge] (a1)--(q)--(a2)--(a3)--(a4)--(a5)--(a1);
\draw[blueedge] (a1)--(p)--(a3)--(a5)--(a2);
\draw[blueedge] (a2) to[out=-45,in=-90] (a4) to[out=90,in=45](a1);
\draw[rededge, dashed] (p) -- (q);
\node[vertexc] (aa2) at (3,0) {$a_2$}; 
\node[vertexc] (aa3) at (4,0) {$a_3$};
\node[vertexc] (aa4) at (5,1) {$a_4$};
\node[vertexc] (aa1) at (3,2) {$a_1$};
\node[vertexc] (aa5) at (4,2) {$a_5$};
\node[vertexc] (p1) at (3.35,1.25) {$p$};
\node[vertexc] (q1) at (3,.6) {$q$};
\draw[blueedge] (q1)--(aa2);
\draw[blueedge](aa3)--(aa4);
\draw[blueedge] (aa5)--(aa1)--(p1)--(aa3)--(aa5)--(aa2);
\draw[blueedge] (aa2) to[out=-45,in=-90] (aa4) to[out=90,in=45](aa1);
\draw[rededge, dashed] (p1) -- (q1);
\node[vertexc] (A2) at (0,-3.5) {$a_2$}; 
\node[vertexc] (A3) at (1,-3.5) {$a_3$};
\node[vertexc] (A4) at (2,-2.5) {$a_4$};
\node[vertexc] (A1) at (0,-1.5) {$a_1$};
\node[vertexc] (A5) at (1,-1.5) {$a_5$};
\node[vertexc] (P) at (.35,-2.25) {$p$};
\node[vertexc] (Q) at (0,-2.9) {$q$};
\draw[blueedge] (A5)--(A1)--(Q)--(A2);
\draw[blueedge] (A3)--(A4);
\draw[blueedge] (P)--(A3)--(A5)--(A2);
\draw[blueedge] (A2) to[out=-45,in=-90] (A4) to[out=90,in=45](A1);
\draw[rededge, dashed] (P) -- (Q);
\draw[->,line width=3pt] (.6,-.5)--(.6,-1);
\node at (3.5,-1) {New main vertex classes:};
\node at (3.5,-1.5) {$\{a_1a_2,a_3\}$, $\{p,a_4,a_5\}$};
\node at (3.5,-2) {(above, $p\ne a_3$)};
\node at (3.5,-2.5) {and};
\node at (3.5,-3) {$\{a_1,a_2,a_3\}$, $\{q,a_4,a_5\}$};
\node at (3.5,-3.5) {(left, $p$ may be $a_3$)};
\node at (3.5,-4.2) {\ \  };
\draw[->,line width=3pt] (2.2,.7)--(2.7,.7);
\end{tikzpicture} 
\quad
\begin{tikzpicture} 
\node[vertexc] (a2) at (0,0) {$a_2$};
\node[vertexc] (a3) at (1,0) {$a_3$};
\node[vertexc] (a4) at (2,1) {$a_4$};
\node[vertexc] (a1) at (0,2) {$a_1$};
\node[vertexc] (a5) at (1,2) {$a_5$};
\node[vertexc] (p) at (0,1) {$p$};
\node[vertexc] (q) at (1,1) {$q$};
\draw[blueedge] (a1)--(p)--(a2)--(a3)--(a4)--(a5)--(a1);
\draw[blueedge] (a1)--(a3)--(q)--(a5)--(a2);
\draw[blueedge] (a2) to[out=-45,in=-90] (a4) to[out=90,in=45](a1);
\draw[rededge, dashed] (p) -- (q); 
\node[vertexc] (aa2) at (3,0) {$a_2$};
\node[vertexc] (aa3) at (4,0) {$a_3$};
\node[vertexc] (aa4) at (5,1) {$a_4$};
\node[vertexc] (aa1) at (3,2) {$a_1$};
\node[vertexc] (aa5) at (4,2) {$a_5$};
\node[vertexc] (p1) at (3,1) {$p$};
\node[vertexc] (q1) at (4,1) {$q$};
\draw[blueedge] (aa1)--(p1)--(aa2)--(aa3)--(aa4)--(aa5);
\draw[blueedge] (aa3)--(q1)--(aa5)--(aa2);
\draw[blueedge](aa4) to[out=90,in=45](aa1);
\draw[rededge, dashed] (p1) -- (q1);
\node[vertexc] (A2) at (0,-3.5) {$a_2$}; 
\node[vertexc] (A3) at (1,-3.5) {$a_3$};
\node[] (A4) at (2,-2) {};
\node[vertexc] (A1) at (0,-1.5) {$a_1$};
\node[vertexc] (A5) at (1,-1.5) {$a_5$};
\node[vertexc] (P) at (0,-2.5) {$p$};
\node[vertexc] (Q) at (1,-2.5) {$q$};
\draw[blueedge] (A5)--(A1)--(P)--(A2)--(A3);
\draw[blueedge] (A1)--(A3)--(Q)--(A5)--(A2);
\draw[rededge, dashed] (P) -- (Q);
\node at (3,-1) {New main vertex classes:};
\node at (3,-1.5) {$\{a_4,q,a_2\}$, $\{a_3,p,a_5\}$};
\node at (3,-2) {(above)};
\node at (3,-2.5) {and};
\node at (3,-3) {$\{a_1,q,a_2\}$, $\{a_3,p,a_5\}$};
\node at (3,-3.5) {(left)};
\node at (3,-4.2) {\ \  };
\draw[->,line width=3pt] (.6,-.5)--(.6,-1);
\draw[->,line width=3pt] (2.2,.7)--(2.7,.7);
\end{tikzpicture} 
\caption{\small Illustration for Lemma~\ref{lm:onecross} when $X$ is a subdivided $K_5$. The resulting new cross-inducing subgraphs are both a
subdivided $K_{3,3}$, the new main vertex classes listed for each. In the first case considered we may have $p=a_3$. In the second case considered we may have $p = a_1$. Otherwise, $p,q$ are subdividing vertices.} 
\label{fig:geza1}
\end{figure}
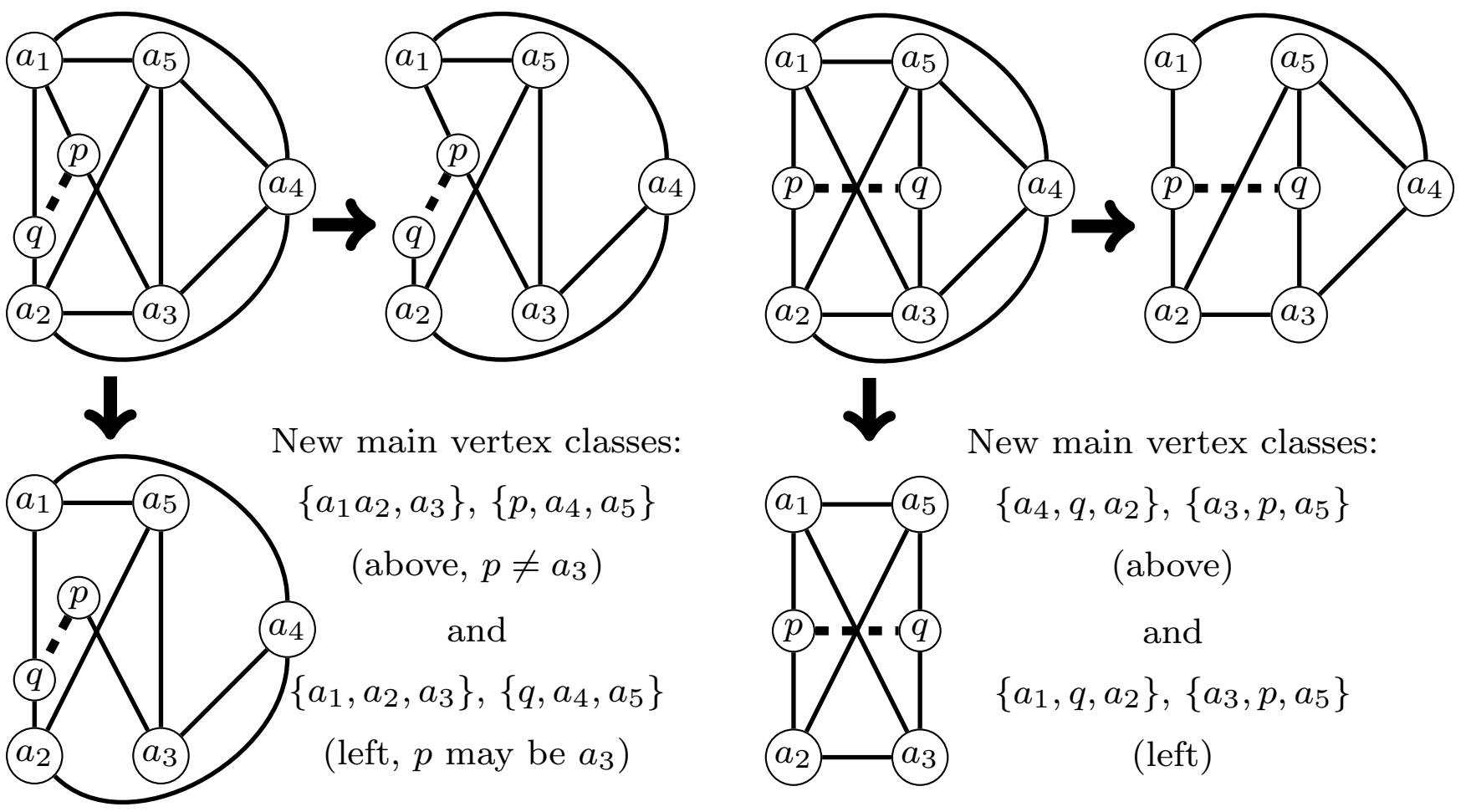

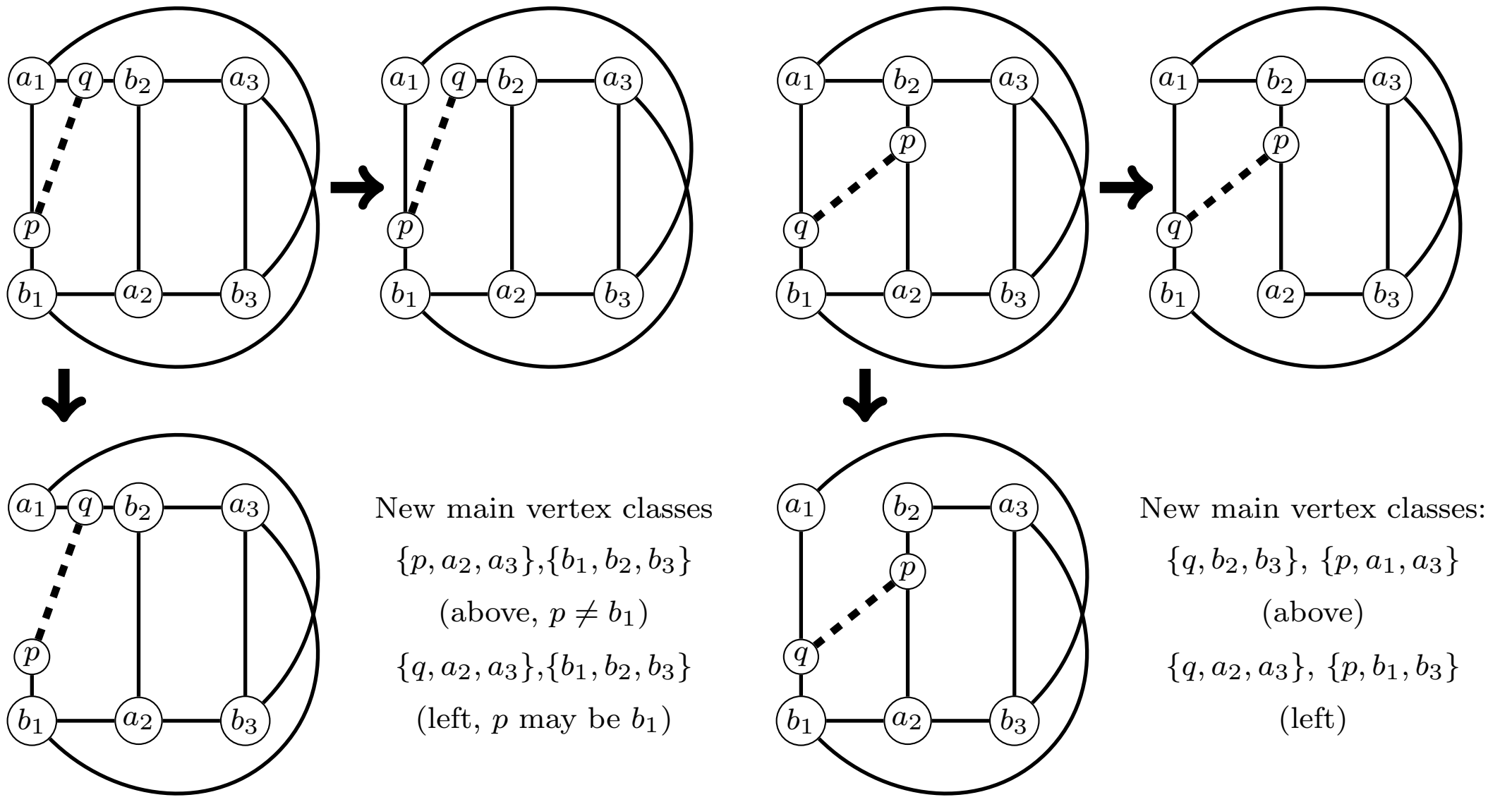
\begin{figure}[http]
\begin{tikzpicture} 
\node[vertexc] (b1) at (0,0) {$b_1$};
\node[vertexc] (b2) at (1,2) {$b_2$};
\node[vertexc] (b3) at (2,0) {$b_3$};
\node[vertexc] (a1) at (0,2) {$a_1$};
\node[vertexc] (a2) at (1,0) {$a_2$};
\node[vertexc] (a3) at (2,2) {$a_3$};
\node[vertexc] (q) at (.5,2) {$q$};
\node[vertexc] (p) at (0,.6) {$p$};
\draw[blueedge] (a2)--(b1)--(p)--(a1)--(q)--(b2)--(a2)--(b3)--(a3)--(b2);
\draw[blueedge] (b1) to[out=-45,in=-135] (2.3,-.3) to[out=45,in=-45] (a3);
\draw[blueedge] (a1) to[out=45,in=135] (2.3,2.3) to[out=-45,in=45] (b3);
\draw[rededge, dashed] (p) -- (q);
\node[vertexc] (bb1) at (3.5,0) {$b_1$}; 
\node[vertexc] (bb2) at (4.5,2) {$b_2$};
\node[vertexc] (bb3) at (5.5,0) {$b_3$};
\node[vertexc] (aa1) at (3.5,2) {$a_1$};
\node[vertexc] (aa2) at (4.5,0) {$a_2$};
\node[vertexc] (aa3) at (5.5,2) {$a_3$};
\node[vertexc] (q1) at (4,2) {$q$};
\node[vertexc] (p1) at (3.5,.6) {$p$};
\draw[blueedge] (aa2)--(bb1)--(p1)--(aa1);
\draw[blueedge] (q1)--(bb2)--(aa2)--(bb3)--(aa3)--(bb2);
\draw[blueedge] (bb1) to[out=-45,in=-135] (5.8,-.3) to[out=45,in=-45] (aa3);
\draw[blueedge] (aa1) to[out=45,in=135] (5.8,2.3) to[out=-45,in=45] (bb3);
\draw[rededge, dashed] (p1) -- (q1);
\draw[->,line width=3pt] (2.8,1)--(3.3,1);
\node[vertexc] (B1) at (0,-4) {$b_1$}; 
\node[vertexc] (B2) at (1,-2) {$b_2$};
\node[vertexc] (B3) at (2,-4) {$b_3$};
\node[vertexc] (A1) at (0,-2) {$a_1$};
\node[vertexc] (A2) at (1,-4) {$a_2$};
\node[vertexc] (A3) at (2,-2) {$a_3$};
\node[vertexc] (Q) at (.5,-2) {$q$};
\node[vertexc] (P) at (0,-3.4) {$p$};
\draw[blueedge] (A2)--(B1)--(P);
\draw[blueedge] (A1)--(Q)--(B2)--(A2)--(B3)--(A3)--(B2);
\draw[blueedge] (B1) to[out=-45,in=-135] (2.3,-4.3) to[out=45,in=-45] (A3);
\draw[blueedge] (A1) to[out=45,in=135] (2.3,-1.7) to[out=-45,in=45] (B3);
\draw[rededge, dashed] (P) -- (Q);
\node at (4.8,-2) {New main vertex classes};
\node at (4.8,-2.5){$\{p,a_2,a_3\}$,$ \{b_1,b_2,b_3\}$};
\node at (4.8,-3){(above, $p\ne b_1$)};
\node at (4.8,-3.5) {$\{q,a_2,a_3\}$,$ \{b_1,b_2,b_3\}$};
\node at (4.8, -4) {(left, $p$ may be $b_1$)};
\draw[->,line width=3pt] (.3,-.7)--(.3,-1.2);
\end{tikzpicture} 
\quad
\begin{tikzpicture} 
\node[vertexc] (b1) at (0,0) {$b_1$};
\node[vertexc] (b2) at (1,2) {$b_2$};
\node[vertexc] (b3) at (2,0) {$b_3$};
\node[vertexc] (a1) at (0,2) {$a_1$};
\node[vertexc] (a2) at (1,0) {$a_2$};
\node[vertexc] (a3) at (2,2) {$a_3$};
\node[vertexc] (p) at (1,1.4) {$p$};
\node[vertexc] (q) at (0,.6) {$q$};
\draw[blueedge] (a2)--(b1)--(q)--(a1)--(b2)--(p)--(a2)--(b3)--(a3);
\draw[blueedge] (b2)--(a3);
\draw[blueedge] (b1) to[out=-45,in=-135] (2.3,-.3) to[out=45,in=-45] (a3);
\draw[blueedge] (a1) to[out=45,in=135] (2.3,2.3) to[out=-45,in=45] (b3);
\draw[rededge, dashed] (p) -- (q);
\node[vertexc] (bb1) at (3.5,0) {$b_1$}; 
\node[vertexc] (bb2) at (4.5,2) {$b_2$};
\node[vertexc] (bb3) at (5.5,0) {$b_3$};
\node[vertexc] (aa1) at (3.5,2) {$a_1$};
\node[vertexc] (aa2) at (4.5,0) {$a_2$};
\node[vertexc] (aa3) at (5.5,2) {$a_3$};
\node[vertexc] (p1) at (4.5,1.4) {$p$};
\node[vertexc] (q1) at (3.5,.6) {$q$};
\draw[blueedge] (bb1)--(q1)--(aa1)--(bb2)--(p1)--(aa2)--(bb3)--(aa3);
\draw[blueedge] (bb2)--(aa3);
\draw[blueedge] (bb1) to[out=-45,in=-135] (5.8,-.3) to[out=45,in=-45] (aa3);
\draw[blueedge] (aa1) to[out=45,in=135] (5.8,2.3) to[out=-45,in=45] (bb3);
\draw[rededge, dashed] (p1) -- (q1);
\node at (4.8,-2) {New main vertex classes:};
\node at (4.8,-2.5) {$\{q,b_2,b_3\}$, $\{p,a_1,a_3\}$};
\node at (4.8,-3) {(above)};
\node at (4.8,-3.5) {$\{q,a_2,a_3\}$, $\{p,b_1,b_3\}$};
\node at (4.8,-4) {(left)};
\draw[->,line width=3pt] (2.8,1)--(3.3,1);
\node[vertexc] (B1) at (0,-4) {$b_1$}; 
\node[vertexc] (B2) at (1,-2) {$b_2$};
\node[vertexc] (B3) at (2,-4) {$b_3$};
\node[vertexc] (A1) at (0,-2) {$a_1$};
\node[vertexc] (A2) at (1,-4) {$a_2$};
\node[vertexc] (A3) at (2,-2) {$a_3$};
\node[vertexc] (P) at (1,-2.6) {$p$};
\node[vertexc] (Q) at (0,-3.4) {$q$};
\draw[blueedge] (A2)--(B1)--(Q)--(A1);
\draw[blueedge] (A3)--(B2)--(P)--(A2)--(B3)--(A3);
\draw[blueedge] (B1) to[out=-45,in=-135] (2.3,-4.3) to[out=45,in=-45] (A3);
\draw[blueedge] (A1) to[out=45,in=135] (2.3,-1.7) to[out=-45,in=45] (B3);
\draw[rededge, dashed] (P) -- (Q);
\draw[->,line width=3pt] (.6,-.7)--(.6,-1.2);
\end{tikzpicture} 
\caption{\small Illustration for Lemma~\ref{lm:onecross} when $X$ is a subdivided $K_{3,3}$. The resulting new cross-inducing subgraphs are both a
subdivided $K_{3,3}$, the new main vertex classes listed for each. In the first case considered we may have $p=b_1$. In the second case considered we may have $q = a_1$ or $p = b_2$. Otherwise $p,q$ are subdividing vertices.}
\label{fig:geza2}
\end{figure}
\end{proof}


\begin{lemma}\label{lm:k331} 
 Let $G$ be a graph, $X$ a subdivided $K_{3,3}$ in $G$ with main vertex classes $\{a_1, a_2, a_3\}$ and $\{b_1, b_2, b_3\}$.
\begin{enumerate}[label={\upshape (\roman*)}]
\item\label{case:starpath} If there is a vertex $v$ in $G-X$ and $3$ internally vertex-disjoint $(G, X)$-paths $Q_1, Q_2, Q_3$ from $v$ to $a_1,a_2,a_3$ respectively, 
then $G$ has at least $4$ cross-inducing subgraphs.
\item\label{case:k5generate} If there are three internally vertex-disjoint $(G,X)$-paths, $P_1$ between $a_1$ and $a_3$, $P_2$
between $a_2$ and $a_3$ and $P_3$ between $b_1$ and $b_2$, then $G$ contains a subdivided $K_5$.
\end{enumerate}
\end{lemma}

\begin{proof}
\ref{case:starpath}: We get a new cross-inducing subgraph for each $i \in [3]$ that is topologically isomorphic to $K_{3,3}$ and has main vertex classes 
$\{a_1, a_2, a_3\}$ and $(\{b_1, b_2, b_3\} \cup \{v\}) \setminus \{b_i\}$ by removing $b_i$ and the internal vertices of the $X$-edge-paths starting from 
$b_i$ and including the paths $Q_1, Q_2, Q_3$.

\ref{case:k5generate}: We get a cross-inducing subgraph topologically isomorphic to $K_5$ on main vertices $a_1, a_2, a_3, b_1, b_2$ by adding to $X$ the paths 
$P_1, P_2$ and $P_3$, and removing internal vertices of the $X$-edge-path
connecting $b_3$ to $a_3$.
\end{proof}


\section{Structure in 2-connected graphs}\label{sec:structure}

In the following lemmas, we add $2$-connectedness of the graph to our conditions. In most  cases, we heavily rely on the structures introduced in Definition~\ref{def:k33structure}.

\begin{lemma}\label{lm:K5linear} 
Let $G$ be a $2$-connected graph that contains a subdivided $K_5$. Then, $G = \X_G$. Furthermore, if $G$ contains exactly $k$ cross-inducing subgraphs, then $G$ has at most $3k+7$ edge-paths.
\end{lemma}

\begin{proof}
 Let $X$ be a subdivided $K_5$ of $G$.
 Let $e$ be any edge of $G$ not in $X$, and let $s,t$ be the endvertices of $e$. Since $G$ is $2$-connected, there are disjoint $(G, X)$-paths $P_1, P_2$ connecting $s, t$ to two different vertices $x, y$ of $X$. Then $P_1eP_2$ is a $(G, X)$-path between two vertices of $X$, and by Lemma~\ref{lm:onecross}, $X\cup P_1eP_2 \subseteq \X_G$. This implies that every edge of $G$ is in $\X_G$, consequently $G \subseteq \X_G \subseteq G$. The first part of the statement follows.

 Since $G$ is 2-connected, we can construct a sequence $X = G_0\subseteq G_1 \subseteq ... \subseteq G_m = G$ where $G_i$ is obtained from $G_{i-1}$ by adding a path $P_i$ between two vertices $x_i$ and $y_i$ of $G_{i-1}$ such that $V(P_i)\cap V(G_{i-1}) = \{x_i,y_i\}$. By 2-connectedness, each $G_{i-1}$ has disjoint paths $P,Q$ from $\{x_i,y_i\}$ to $X$. Thus, $G_{i}$ has a cross-inducing subgraph using all edges of $P_i$, and has at least one more cross-inducing subgraph than $G_{i-1}$. Furthermore, $G_i$ has at most $3$ more edge-paths than $G_{i-1}$. 
 So, $G_m$ has at most $3m + 10$ edge-paths. Further, $G_m$ has at least $m+1$ cross-inducing subgraphs, but since $G_m = G$ it must be that $m+1 \leq k$. So $3m + 10 \leq 3(k-1) + 10 = 3k+7$.
\end{proof}

\begin{lemma}\label{lm:edgeestimate0}Let G be a 2-connected graph that contains no subdivided $K_5$ and let $X$ be a cross-inducing subgraph of $G$.
\begin{enumerate}[label={\upshape (\roman*)}]
\item\label{case:decompG} $G = G(X) \cup E^*_X \cup \Big( \bigcup_{H\in \mathcal{P}_1^X }\hat{H}\Big)$
\item\label{case:outsideEdges}For each $e\in F^*_X$, the graph $X\cup\{e\}$ is 2-connected and has a cross-inducing subgraph that contains the edge $e$.
\end{enumerate}

\end{lemma}
\begin{proof} \ref{case:decompG} is obvious from the definitions.

\ref{case:outsideEdges} Each $e\in F_x^*$ is a $(G,X)$-path between two vertices of $X$ that are not $X$-cooperative, resulting in a new cross-inducing subgraph containing $e$ in $X\cup \{e\}$ by Lemma \ref{lm:onecross}. As $X$ is 2-connected and for any edge $e\in F_X^*$, $e$ is in a cycle in $X\cup\{e\}$, the statement follows.
\end{proof}

\begin{lemma}\label{lm:bipartiteTruths} Let $G$ be a 2-connected graph that contains no subdivided $K_5$ and let $X$ be a cross-inducing subgraph of $G$. Let $H\in\mathbb{P}_X$. The following are true.
\begin{enumerate}[label={\upshape (\roman*)}]
\item\label{case:decomposition} $|N_H|\geq 2$. Moreover, for any fixed $h\in N_H$, every edge of $\hat{H}$ is on some $(\hat{H},X)$-path connecting $h$ to some vertex $s\in N_H-\{h\}$ and $X\cup\hat{H}$ is 2-connected.
\item\label{case:triplePaths} If $|N_H|\geq 3$, then for any two vertices $h_1,h_2\in N_H$ and edge $e$ of $\hat{H}$, there is a vertex $z\in H$, a vertex $h_3\in N_H-\{h_1,h_2\}$, and internally disjoint $(\hat{H},X)$-paths $P_1,P_2$, and $P_3$ from $z$ to $h_1,h_2,$ and $h_3$, respectively such that one of the $P_i$ contains $e$. Moreover, for any $3$ vertices $s_1,s_2,s_3\in N_H$, there is a vertex $z^*$ and internally disjoint paths $R_1,R_2$, and $R_3$ from $z$ to $s_1,s_2,$ and $s_3$, respectively.
\item\label{case:subdividingNeighbor} If $N_H$ contains a subdividing vertex $h$ of $X$, then every edge of $\hat{H}$ is on a $(\hat{H},X)$-path $P$ connecting $h$ to some vertex of $X$, and every edge of $\hat{H}$ is on a cross-inducing subgraph of $X\cup \hat{H}$.
\item\label{case:noncooperativeNeighbors} If $N_H$ contains only $X$-main vertices and $N_H$ contains two vertices that are not $X$-cooperative, then every edge of $\hat{H}$ is on a $(\hat{H},X)$-path $P$ between two vertices of $N_H$ that are not $X$-cooperative, so every edge of $\hat{H}$ is on a cross-inducing subgraph of $X\cup\hat{H}$.
\item\label{case:mainClassNeighbors} If $N_H = \{t_1,t_2,t_3\}$ where $\{t_1,t_2,t_3\}$ is a main vertex class of $X$, then for every edge of $\hat{H}$ there is a vertex $z$ of $H$ and internally disjoint $(\hat{H},X)$ paths $P_1,P_2$, and $P_3$ from $z$ to $t_1,t_2,$ and $t_3$, respectively, such that $e$ is on one of the $P_i$. Consequently, every edge of $\hat{H}$ is on a cross-inducing subgraph of $X\cup\hat{H}$.
\item\label{case:containment} $G(X)$ is 2-connected and $G(X)= \mathcal{X}_{G(X)} \subseteq \mathcal{X}_G$
\end{enumerate}
\end{lemma}
\begin{proof}
\ref{case:decomposition}: As $H\cup F_H \subseteq \hat{H}$, $\hat{H}$ is a connected graph that is not edgeless. First we will show the weaker claim that every edge of $\hat{H}$ is on a $(\hat{H},X)$-path connecting two vertices of $X$. This immediately gives $|N_H|\geq 2$ and, as $X$ is 2-connected, so is $X\cup\hat{H}$ .

Let $e$ be an edge of $\hat{H}$ with endpoints $x,y$. Since $G$ is 2-connected, there are disjoint $(G,X)$-paths $P$ and $Q$ such that $P$ connects $x$ to some $s_1\in X$ and $Q$ connects $y$ to some $s_2\in X$, $s_1\ne s_2$. As $H$ is a component of $G-X$, the internal vertices of $P,Q$ are in $H$. $PeQ$ is a $(\hat{H},X)$-path between two vertices of $N_H$ that contains the edge $e$. The weaker claim follows.

Now fix an $h\in N_H$. Consider an edge $e$ of $\hat{H}$. By our weaker claim, $e$ lies on a $(\hat{H},X)$-path $P$ between two vertices $s_1,s_2 \in N_H$; then $P$ has at least one internal vertex. If $h \in \{s_1,s_2\}$, we are done. Otherwise let $P^0$ be the path $P$ without the endvertices $s_1,s_2$. Then $P^0$ is a path in $H$ with at least one vertex. Let $f\in F_H$ be an edge connecting $h$ to some $v \in H$. As $H$ is connected, there is an $(H,P)$-path $Q$ from $v$ to $P^0$ ending on a vertex $z$. Now $fQzPs_1$, $fQzPs_2$ are both $(\hat{H},X)$-paths between $h$ and a vertex of $N_H - \{h\}$, and one of these paths contain $e$. The claim follows.

\ref{case:triplePaths}: Fix $h_1,h_2$, and let $e$ be an edge of $\hat{H}$. By \ref{case:decomposition} there is a $(\hat{H},X)$-path $Q$ from $h_1$ to some $s \in N_H$ containing the edge $e$. If $s\ne h_2$, set $t = h_2$, otherwise choose a $t \in N_H -\{h_1,h_2\}$ arbitrarily. Let $f \in F_H$ 
be an edge with endvertices $t, y$; $y \in H$. As $H$ is connected, there is an $(H,Q)$-path $R$ from $y$ to a vertex $z$ of $Q \in H$. The required paths $P_1,P_2,P_3$ are $h_1Qz$,$zQs$ and $Rf$ (not necessarily in the same order).

Now let $s_1,s_2,s_3\in N_H$. We follow the proof of the first part of the claim: it gives that there is a $(\hat{H},X$)-path $L$ connecting $s_1,s_2$,
and we have a vertex $z^*$ and internally disjoint $(\hat{H},X)$-paths $R_1,R_2,R_3$ from $z^*$ to each of $s_1,s_2,s_3$.

\ref{case:subdividingNeighbor}: By ~\ref{case:decomposition}, every edge $e$ of $\hat{H}$ is on a $(\hat{H},X)$-path between $h$ and some vertex $s \in N_H - \{h\}$. Lemma~\ref{lm:onecross} yields that $e$ is on a cross-inducing subgraph of $X \cup P$.

\ref{case:noncooperativeNeighbors}: Let $t_1,t_2$ be two vertices of $N_H$ that are not $X$-cooperative. Let $e$ be an edge in $\hat{H}$. If $|N_H| = 2$, then by \ref{case:decomposition}, $e$ is on a $(\hat{H},X)$-path $P$ connecting $t_1$ and $t_2$, which are not $X$-cooperative. By Lemma \ref{lm:onecross} we get that $e$ is on a cross-inducing subgraph of $X\cup P$. If $|N_H|\geq 3$ then by \ref{case:triplePaths} there is a vertex $v_e$ of $H$, a vertex $t_3 \in N_H -\{t_1,t_2\}$, and internally vertex disjoint paths $Q_1, Q_2, Q_3$ from $v_e$ to $t_1, t_2, t_3$ respectively such that $e$ is one one of the $Q_i$. Without loss of generality, $t_2$ and $t_3$ are $X$-cooperative. Then $Q_1Q_2$, $Q_1Q_3$ are $(\hat{H},X)$-paths connecting two $X$-main vertices that are not $X$-cooperative, and at least one of these paths (call it $P$) contains $e$. By Lemma \ref{lm:onecross} we get that $e$ is on a cross-inducing subgraph of $X \cup P$.

\ref{case:mainClassNeighbors} Let $e$ be an edge of $\hat{H}$. By \ref{case:triplePaths} there is a vertex $v_e$ and internally disjoint paths $P_1,P_2,P_3$ connecting $v_e$ to the three vertices of $N_H$ such that $e$ is on one of these paths. By Lemma \ref{lm:k331} \ref{case:starpath}, $e$ is on a cross-inducing subgraph of $X \cup (P_1P_2) \cup P_3$.

\ref{case:containment}: Let $D$ be the 2-connected component of $G(X)$ containing $X$. By \ref{case:decomposition} and Lemma \ref{lm:edgeestimate0} \ref{case:outsideEdges}, $G(X) \subseteq D \subseteq G(X)$, so $G(X)$ is 2-connected.

Let $H \in \mathcal{P}_2^X$ , then one of the following must hold:
\begin{enumerate}
    \item[(a)] $N_H$ contains a subdividing vertex of $X$
    \item[(b)] $N_H$ contains two $X$-main vertices $s,t$ that are not $X$-cooperative
    \item[(c)] $N_H = \{s, h, t\}$, where $s$, $h$ and $t$ are pairwise $X$-cooperative.
\end{enumerate}

By \ref{case:subdividingNeighbor}, \ref{case:noncooperativeNeighbors}, \ref{case:mainClassNeighbors} and Lemma \ref{lm:edgeestimate0} \ref{case:outsideEdges}, we get $G(X) \subseteq \mathcal{X}_{G(X)} \subseteq G(X)$. The statement follows.
\end{proof}

\begin{lemma}\label{lm:k33Metric}
Let $G$ be a 2-connected graph with exactly $k$ cross-inducing subgraphs, none of which are a subdivided $K_5$. Let $X$ be a cross-inducing subgraph of $G$. Set $k_0$ be the number of cross-inducing subgraphs in $G(X)$ and for each $H \in \mathcal{P}_1^X$, let $k_H$ be the number of cross-inducing subgraphs in $\hat{H}\cup \gamma_{s_H,t_H}$. We have the following:

\begin{enumerate}[label={\upshape (\roman*)}]
\item\label{case:connectivityOfUnion} For each $H \in \mathcal{P}_1^X$ the graph $\hat{H}\cup \gamma_{s_H,t_H}$ is 2-connected, and every edge in $G-H$ is on a $(G,\hat{H})$-path connecting $s_H$ and $t_H$.
\item\label{case:kInequality} $1\leq k_0 \leq k$
and $k_0 + \sum\limits_{H\in \mathcal{P}_1^X} k_H \leq k$
\item\label{case:measureInequality} If for some $H\in\mathcal{P}_1^X$ the graph $\hat{H}\cup \gamma_{s_H,t_H}$ has a cross-inducing subgraph that contains the edge $\gamma_{s_H,t_H}$, then $G$ has a cross-inducing subgraph $X^*$ with $m(X^*)> m(X)$.
\end{enumerate}
\end{lemma}
\begin{proof}

\ref{case:connectivityOfUnion}: Let $H\in\mathcal{P}_1^X$. $|\hat{H}|=|H|+2\geq 3$. 
By Lemma \ref{lm:bipartiteTruths} \ref{case:decomposition} every edge $f$ of $\hat{H}$ is on a path $P$ that connects $s_H$ and $t_H$ in $\hat{H}$; therefore $P\cup \gamma_{s_H, t_H}$ is a cycle containing $f$ and $\gamma_{s_H, t_H}$ . This implies that $\hat{H}\cup \gamma_{s_H,t_H}$ is 2-connected.

Let $e$ be an edge of $G - H$, and let $x,y$ be the endvertices of $e$. As $G$ is 2-connected, there are disjoint paths $P,Q$ connecting $\{x,y\}$ to $\{s_H,t_H\}$. As ${s_H,t_H}$ is a cutset, we must have that $P,Q$ are paths in $G - H$, and $P \cup {e} \cup Q$ is a $(G,\hat{H})$-path containing $e$. We are done.

\ref{case:kInequality}: $1\leq k_0 \leq k$ follows from $X\subseteq G(X)$ and Lemma \ref{lm:bipartiteTruths} \ref{case:containment}. For each $H\in\mathcal{P}_1^X$, let $P_H$ be a path in $X$ connecting $s_H$ and $t_H$, $\mathcal{Z}_H$ be the set of cross-inducing subgraphs of $\hat{H}\cup \gamma_{s_H,t_H}$ that contain the edge $\gamma_{s_H,t_H}$, and set $\mathcal{Z}_H^* = \{(Z-\gamma_{s_H,t_H})\cup P_H \ : \ Z\in \mathcal{Z}_H\}$. Setting $k_H^*$ as the number of cross-inducing subgraphs of $\hat{H}$, we get that $k_H = k_H^* + |\mathcal{Z}_H|$. Clearly, $|\mathcal{Z}_H| = |\mathcal{Z}_H^*|$, $\mathcal{Z}_H^*$ is a collection of cross-inducing subgraphs of $G$, and these subgraphs are not in 
any $\hat{H'}\cup \gamma_{s_{H'},t_{H'}}$ for $H'\ne H$.

For $H'\ne H$ we have $\mathcal{Z}_H^* \cap \mathcal{Z}_{H'}^* = \emptyset$. This implies that
$$k_0 + \sum\limits_{H\in \mathcal{P}_1^X} k_H = k_0 + \sum\limits_{H\in \mathcal{P}_1^X} (k_H^* + |\mathcal{Z}_H|) = k_0 + \sum\limits_{H\in \mathcal{P}_1^X} (k_H^* + |\mathcal{Z}_H^*|) \leq k$$

and the claim follows.

\ref{case:measureInequality}: Assume for $H \in \mathcal{P}_1^X$ the graph $\hat{H}\cup \gamma_{s_H,t_H}$  contains a cross-inducing subgraph $X'$ that contains the edge $\gamma_{s_H,t_H}$. Clearly, $m(X) \leq |G| - |X| - |H|$, and the vertices $s_H,t_H$ lie on the same edge-path of $X'$. Let $P$ be an $s_H$-$t_H$ path in $X \cup E_X^*$ that contains the fewest $X$-main vertices. (Note that such a path exists in $X$ and covers at most three of the six $X$-main vertices, consequently $|P|<|X|$). Let $X^* = (X' - \gamma_{s_H,t_H}) \cup P$. Then $X^*$ is a subdivision of $X'$, so $X^*$ is a cross-inducing subgraph of $G$, consequently both $X'$ and $X^*$ are subdivided $K_{3,3}$'s. Let $g$ be an edge of $G - H$ such that $g$ is not an edge of $P$. If both endvertices of $g$ are on $P$, then $g$ is a $(G, H \cup P )$-path between two vertices $x,y$ of $P$. 
Otherwise by part \ref{case:connectivityOfUnion} $g$ is on a $(G,\hat{H})$-path $Q$ connecting $s_H$ and $t_H$. This implies
that (by sufficiently shortening $Q$ if needed) $g$ is on a $(G, \hat{H} \cup P )$-path between two vertices $x,y$ of $P$. As $x,y$ are on the same $X^*$-edge-path, $g$ is an edge of a $\widehat{H^*}$ for some $H^* \in \mathcal{P}_2^{X^*}$. Consequently $m(X^*) \geq |G| - |\hat{H}\cup P| > |G| - |H\cup X| \geq m(X)$.
\end{proof}


Now we are ready to prove the following
\begin{lemma}\label{lm:meta}
Let $g:\mathbb{N}\rightarrow\mathbb{N}$ be a superadditive function and $k\in\mathbb{N}$. We have that $f(k)\le g(k)$
if the following holds all $2$-connected graphs $G$ 
with at most $k$ cross-inducing subgraphs and all $1\le r\le k$
\begin{enumerate}[label={\upshape(\arabic*)}]
    \item\label{case:k5norestrict} If $G$  contains a subdivided $K_5$, and has exactly $r$ cross-inducing subgraphs, $G$ has a drawing with at most $g(r)$ crossings
    \item\label{case:reducedtok33}  
    If $G$ contains no subdivided $K_5$, and $G=G(X)\cup S$ for some cross-inducing subgraph $X$ and set of edges $S$ that connect $X$-cooperative vertices, and $G(X)$ has exactly $r$ cross-inducing subgraphs, then $G$ has a drawing with at most $g(r)$ crossings where only edges of $G(X)$ cross.
\end{enumerate}
\end{lemma}

Note that while our definition of graphs allow multiple edges, the set $S$ in part~\ref{case:reducedtok33} is a set and not a multiset, so there is at most one edge between any two $X$-cooperative vertices.
We also observe that for  $r=k$  condition~\ref{case:reducedtok33} implies that 
no edges of $S$ lie in a cross-inducing subgraph of $G$.

\begin{proof}
By Lemma~\ref{lm:noYcross}, for $f(k)\le g(k)$ it is enough to show that any $2$-connected graph $G$ with $k^*\le k$ cross-inducing subgraphs has a drawing with at most $g(k^*)$ crossings in which only edges of $\X_G$ cross. 

If $k^*=0$, Kuratowski Theorem implies that the graph has a drawing with $0$ crossings, and as $0\le g(0)$, the statement holds. We may assume that $k^*\ge 1$.

If $G$ contains a subdivided $K_5$,  Lemma~\ref{lm:K5linear} gives that $\X_G=G$, and we are done by~\ref{case:k5norestrict}. 

If $G$ does not contain a subdivided $K_5$, we will use induction on $k^*$. 

Assume that for every $k'<k^*$ every $2$-connected graph $G^*$ that does not contain a $K_5$ and has  exactly $k'$
cross-inducing subgraphs has a drawing with at most $g(k')$ crossings in which only edges of $\X_{G^*}$ cross.

Let $G$ be a $2$-connected graph without a $K_5$ that has exactly $k^*$ cross-inducing subgraphs. Let $X$ be a cross-inducing subgraph of $G$ with $m(X)$ maximum among all cross-inducing subgraphs. Set $k_0$ be the number of cross-inducing subgraphs of $G(X)$, and for each $H\in\mathcal{P}_1^X$ let $k_H$ be the number of cross inducing subgraphs of $W_H := \hat{H}\cup\gamma_{s_H,t_h}$. 
By Lemma \ref{lm:bipartiteTruths} \ref{case:decomposition} and \ref{case:containment}, $G(X)$ is 2-connected, $G(X) = \mathcal{X}_{G(X)} \subseteq \mathcal{X}_G$, and each $W_H$ is 2-connected. By Lemma \ref{lm:k33Metric} \ref{case:kInequality} we have $k_0 + \sum\limits_{H\in\mathcal{P}_1^X}k_H \leq k^*$, $1\leq k_0\leq k^*$, for all $H\in\mathcal{P}_2^X$ no cross-inducing subgraphs of $W_H$ contains $\gamma_{s_H,t_h}$, and $k_H < k^*$. 
By the induction hypothesis, each $W_H$ has a drawing
$\mathcal{D}_H$ with at most $g(k_H)$ crossings in which only edges of $\mathcal{X}_{\hat{H}}$ cross. In particular,
$\gamma_{s_H,t_h}$ is not crossed in $\mathcal{D}_H$.

Let $S=\{\gamma_{s,t}:\exists H\in\mathcal{P}_1^X\,\, N_H=\{s,t\}\hbox{ or }\exists e\in E^*_X\hbox{ connecting $s$ and }t\}$.
Since every $H\in\mathcal{P}_1^X$ contains a path between the two vertices of $N_H$,
 $G$ contains a subdivided copy of $G(X)\cup S$ and does not contain a subdivided $K_5$, $G(X)\cup S$ has no subdivided $K_5$ and $G(X)\cup S$ has at most $k$ cross-inducing subgraphs.
By ~\ref{case:reducedtok33}, $G(X) \cup S$ has a drawing $\mathcal{D}_0$ with at most $g(k_0)$ crossings in which only edges of $G(X)$ are crossed.

Lemma \ref{lm:obvious} then gives that we have a drawing of $G$ such that only edges of $\mathcal{X}_G$ cross and the drawing has at most $g(k_0)+\sum_H g(k_H)\le g(k^*)$ crossings. Our proof is finished.
\end{proof}

\section{Counting edge-paths}\label{sec:k33Linear}

This section is dedicated to the  proof the following:
\begin{lemma}\label{lm:k33Linear}
    Let $G$ be a $2$-connected graph with exactly $k$ cross-inducing subgraphs, none of which are a subdivided $K_5$, and $X$ be a cross-inducing subgraph of $G$. Then $G(X)$ has at most $3k+6$ edge-paths.
\end{lemma}

The proof of this lemma will easily follow from Lemmas~\ref{lm:k33LinearwithSubdivision},~\ref{lm:k33LinearwithNoncooperative}
and~\ref{lm:k33LinearWithCooperative}. Lemma~\ref{lm:forkingLemma} is a technical lemma for the proof of Lemma~\ref{lm:k33LinearWithCooperative}.

\begin{lemma}\label{lm:k33LinearwithSubdivision}
    Suppose $G$ is a $2$-connected graph, containing no $K_5$ subdivisions, $X$ be a cross-inducing subgraph and $H\in\mathcal{P}_2^X$ such that $N_H$ contains an $X$-subdividing vertex $b$. Then for any subgraph $Q$ of $G$ with $X\subseteq Q$ and $V(Q)\cap V(H) =\emptyset$, if $Q\cup \hat{H}$ 
    has exactly $k'$ more cross-inducing subgraphs than $Q$, then $Q\cup \hat{H}$ has at most $3k'$ more edge-paths than $Q$.
\end{lemma}
\begin{proof}
    
    Let $G, H, Q,$ and $b$ be as described in the claim. 
    We proceed by constructing a sequence of graphs $\emptyset = H_0 \subseteq H_1 \subseteq H_2 \subseteq ... \subseteq H_m = \hat{H}$ which obeys the following conditions for all $i\in [m]$:

    \begin{enumerate}[label={\upshape (\arabic*)}]
        \item\label{case1:Hbasics} $b\in H_i$, $|H_i\cap N_H|\ge 2$, 
        $H\cap H_i\ne\emptyset$, and $H_i$ is connected,
        \item\label{case1:threemore} There are at most three more 
        edge-paths  in 
        $Q\cup H_i$
        than 
        in
        $Q\cup H_{i-1}$,  
        \item\label{case1:onemore} 
        $Q\cup H_i$  has strictly more cross-inducing subgraphs than $Q\cup H_{i-1}$,
        \item\label{case1:eachx} We set $H_i^*:=H_i-(N_H-\{b\})$. For any vertex $x\in H_i^*$, there exists a connected subgraph $H_i^x$ of $H_i^*$ containing both $x$ and $b$ such that every component of $H_i - H_i^x$ contains a vertex of $N_H - \{b\}$,
        \item\label{case1:eachxy} For any two $x,y\in H_i^*$, $x\notin H_i^y$ or $y\notin H_i^x$.
    \end{enumerate}
    
    Using Lemma \ref{lm:bipartiteTruths} \ref{case:subdividingNeighbor}, there exists a $(\hat{H},X)$-path $P$ from $b$ to some vertex of $N_H-\{b\}$ which contains some edge $e_1$ of $\hat{H}$. Let $H_1 = P$, and~\ref{case1:Hbasics} is clearly satisfied. The addition of $P$ to $Q$ can add at most three new edge-paths, so~\ref{case1:threemore} holds. 
    By Lemma \ref{lm:onecross}, $X\cup H_1$ has a cross-inducing subgraph which includes the path $P$. This cross-inducing subgraph contains the edge $e_1$ which does not belong to $Q$, so $Q\cup H_1$ has at least one more cross-inducing subgraph than $Q$. So \ref{case1:onemore} holds.
    For any $x\in H_1^*$, we define $H_1^x$ as the path $bPx$, which easily satisfies~\ref{case1:eachx} and~\ref{case1:eachxy}. 

    Suppose for some $i\ge 1$ we have constructed $H_{i}$ as desired, and $H_{i}\ne \hat{H}$. 

    We will first show that there is a vertex $z_1\in H_{i}\cap H$ and an edge $e_{i+1}$ incident upon $z_1$ such that $e_{i+1}\in E(\hat{H})-E(H_{i})$. (We will denote the other endvertex of $e_1$ by $v$.)

    If $H\not\subseteq H_{i}$,  $E(H)\not\subseteq E(H_{i})$.  
    Since $H$ is connected, $H\cap H_{i}\ne\emptyset$ and $E(H)-E(H_{i})\ne\emptyset$, there must be some vertex $z_1\in H_{i}\cap H\subseteq H_{i}^*=H_{i}-(N_H-\{b\})$ that is incident upon an edge $e_{i+1}\in E(H)-E(H_{i})\subseteq E(\hat{H})-E(H_{i})$. 
    
    Otherwise, $H\subseteq H_i$ and every edge in $E(\hat{H})-E(H_{i})$ has one endpoint in $N_H$ and the other in $H$. In this case, we let $e_{i+1}$ be any edge in $E(\hat{H})-E(H_{i})$ and $z_1$ be its endpoint in $H$. 

    We achieved our first goal.
 
    Next we will find a $(\hat{H},H_{i}\cup X$)-path $R$ connecting $z_1\in H_{i}^*$ to a vertex $z_2\in H_{i}\cup X$ that contains $e_{i+1}$. If $v\in H_{i}\cup X$, simply let $R=e_{i+1}$ and $z_2=v$. If $v\notin H_{i}\cup X$,  by Lemma \ref{lm:bipartiteTruths} \ref{case:decomposition}
    we have $X\cup\hat{H}$ is $2$-connected, consequently there exists a pair of internally disjoint $(\hat{H}\cup X,H_{i}\cup X)$ paths from $v$ to distinct vertices of $H_{i}\cup X$. As at most one of these paths contain the vertex $z_1$, we let $P$ be the path connecting $v$  to $z_1$ that does not contain $e_{i+1}$; then $P$ is a $(\hat{H}\cup X,H_{i}\cup X)$-path not containing the vertex $z_1$ (and consequently not containing the edge $e_{i+1}$). 
    As $N_H$ separates $X$ from $\hat{H}$, $N_H\subseteq X$ and $z_1\notin X$, $P$ is a $(\hat{H},H_{i}\cup X)$-path.
    Set $R=Pe_{i+1}$ and $z_2$ to be the endpoint of $P$ different from $v$.

    We achieved our second goal.
    
    We let $H_{i+1}=H_{i}\cup R$. Then~\ref{case1:Hbasics} and ~\ref{case1:threemore} hold for $H_{i+1}$, and we must now show that $H_{i+1}$ satisfies conditions \ref{case1:onemore}-\ref{case1:eachxy}. Note that either $z_2\in N_H - \{b\}$ or $z_2\in H_{i}^*$. 

     \textit{Case 1:  
     $z_2\in N_H-\{b\}$:}
     
    Since \ref{case1:eachx} applies to $H_{i}$, $H_{i}^{z_1}$ contains an $b-z_1$ path $T$. Note that $T\subseteq H_{i}^*$ and $V(R)\cap V(H_{i}^*)=\{z_1\}$, therefore $bTz_1Rz_2$ is a $(\hat{H},X)$ path from $b$ to $z_2$ containing the edge $e_{i+1}$. By Lemma \ref{lm:onecross}, $(bTz_1Rz_2)\cup X$ contains a cross-inducing subgraph which contains $bTz_1Rz_2$ and thus contains the edge $e_{i+1}$.
    As $(bTz_1Rz_2)\cup X\subseteq Q\cup H_{i+1}$ and $Q\cup H_{i}$ does not contain the edge $e_{i+1}$, \ref{case1:onemore} is satisfied.
    
    Now, for every $x\in H_{i+1}^*$ we define $H_{i+1}^x$ and show that these choices satisfy \ref{case1:eachx} and \ref{case1:eachxy}.
 
    Either $x\in H_{i}^*$ or $x$ is an internal vertex of $R$. If $x\in H_{i}^*$, we set $H_{i+1}^x = H_{i}^x$. Since $H_{i}^x$ satisfies \ref{case2:eachx}, $H_{i+1}^x$ is connected and contains $x$ and $b$. Furthermore, any component of $H_{i+1} - H_{i+1}^x$ contains $z_2$ or is a component of $H_{i} - H_{i}^x$, so \ref{case1:eachx} is satisfied.
    
    Next, if $x\in H_{i+1}^*-H_i^*=H_{i+1}-H_i$, then $x$ is an internal vertex of $R$. We set $H_{i+1}^x$ to be $xRz_1\cup H_i^{z_1}$. This is clearly connected and contains both $x$ and $b$. Furthermore, any component of $H_{i+1} - H_{i+1}^x$ is either a component of $H_{i} - H_{i}^{z_1}$ or contains $z_2$, so \ref{case1:eachx} is satisfied. 
    
    Condition \ref{case1:eachxy} follows easily, so we are done with the case where $z_2\in N_H-\{b\}$.
    
   \textit{Case 2: 
   $z_2\in H_i^*$:}
   
    As \ref{case1:eachx} and \ref{case1:eachxy} hold for $H_i$, we have subgraphs $H_i^{z_1}$ and $H_i^{z_2}$ of $H_i^*$ such that for some $j\in[2]$, $z_{3-j}\notin H_i^{z_{j}}$. By \ref{case1:eachx}, $H_i^{z_j}$ contains an $b-z_j$ path $T$, and every component of $H_i - H_i^{z_j}$ contains some vertex of $N_H-\{b\}$. 
    Let $S$ be a path in $H_i - H_i^{z_j}$ from $z_{3-j}$ to some $y\in N_H-\{b\}$. Now, set $P^* = bTz_jRz_{3-j}Sy$. This $P^*$ is a $(\hat{H},X)$ path from $b$ to $y$ containing the edge $e_{i+1}$. By Lemma \ref{lm:onecross}, $Q\cup H_{i}$ contains a cross-inducing subgraph which contains $P^*$ and thus contains the edge $e_{i+1}$, so \ref{case1:onemore} is satisfied.

    Now, for every $x\in H_{i+1}^*$ we define $H_{i+1}^x$ and show these choices satisfy \ref{case1:eachx} and \ref{case1:eachxy}.
    
    First, suppose $x\in H_i^*$. If $z_1,z_2\in H_i^x$, then we set $H_{i+1}^x = H_i^x \cup R$. Otherwise, we set $H_{i+1}^x = H_i^x$. In both cases $H_{i+1}^x$ is connected and contains $x$ and $b$. 
    As $H_{i+1}=H_i\cup R$, $H_{i+1}^x\in\{H_i^x,H_i^x\cup R\}$, and if $H_{i+1}^x\ne H_i\cup R$ then at least one of $z_1,z_2$ is not in $H_i^x$,
    every component of $H_{i+1} - H_{i+1}^x$ contains component of $H_{i} - H_{i}^x$, and thus contains a vertex of $N_H - \{b\}$.

    Next, if $x\in H_{i+1}^* - H_i^*=H_{i+1}-H_i$, then $x$ is an internal vertex of the path $R$. Recall that $j$ was chosen such that $z_{3-j}\notin H_i^{z_j}$. In this case, we set $H_{i+1}^x = H_i^{z_j} \cup xRz_j$. Clearly, $H_{i+1}^x$ is connected and contains $x$ and $b$. Every component of $H_{i+1} - H_{i+1}^x$ is either a component of $H_i - H_i^{z_j}$, and thus contains a vertex in $N_H-\{b\}$, 
    or contains the component of $H_i - H_i^{z_j}$ which contains $z_{3-j}$, and thus contains a vertex in $N_H-\{b\}$. So \ref{case1:eachx} is satisfied for our choices of $H_{i+1}^x$.

    Finally we must show that \ref{case1:eachxy} holds. Let $x,y\in H_{i+1}^*$ be distinct. If $x,y\in H_{i}^*$, then condition \ref{case1:eachxy} is satisfied since the vertices of $H_{i+1}^x$ and $H_{i+1}^y$ differ from those of $H_{i}^x$ and $H_i^y$ only possibly by internal vertices of $R$. Furthermore, if $x,y$ are both internal vertices of $R$, then if $x\in H_{i+1}^y$ it must be that $x\in z_jRy$. Therefore, $y\notin H_{i+1}^x$. So it suffices to assume that $x\in H_i^*$ and $y$ is an internal vertex of $R$. Suppose $x\in H_{i+1}^y$. Then by construction, $x\in H_i^{z_j}$. Since $x\in H_i^*$, $H_{i+1}^x$ contains $y$ only if $H_i^x$ contains both $z_1$ and $z_2$. But since $x\in H_i^{z_j}$, $z_j\notin H_i^x$. Therefore $y\notin H_{i+1}^x$, and we are done.

    Therefore, we have the desired sequence. Note that by~\ref{case1:onemore}, $Q\cup H_m = Q\cup \hat{H}$ contains at least $m$ more cross-inducing subgraphs than $Q$, and by~\ref{case1:threemore}, $Q\cup H_m$ has at most $3m$ more edge-paths than $Q$. Since $Q \cup H_m$ has $k'$ more cross inducing subgraphs than $Q$, $m \leq k'$. Therefore, $Q \cup H_m$ has at most $3m \leq 3k'$ more edge-paths than $Q$. This completes the proof.

\end{proof}

\begin{lemma}\label{lm:k33LinearwithNoncooperative}
    Suppose $G$ is a 2-connected graph, containing no $K_5$ subdivisions, $X$ be a cross-inducing subgraph of $G$, $H\in\mathcal{P}_2^X$ such that $N_H$ contains a pair of $X$-main vertices that are not $X$-cooperative, and $N_H$ contains no $X$-subdividing vertices. Then for any subgraph $Q$ of $G$ with $X\subseteq Q$ and $V(Q)\cap V(H) =\emptyset$, if $Q\cup \hat{H}$ 
    has exactly $k'$ more cross-inducing subgraphs than $Q$, then $Q\cup \hat{H}$ has at most $3k'$ more edge-paths than $Q$.
\end{lemma}
\begin{proof}
    
    Let $G, H,$ and $Q$ be as described in the claim. Denote the main vertex classes of $X$ by $A$ and $B$. By assumption, $N_H\subseteq A\cup B$, and neither $N_H\cap A$ nor $N_H\cap B$ is empty.
    
    Similarly to the previous proof, we proceed by constructing a sequence of graphs $\emptyset = H_0 \subseteq H_1 \subseteq H_2 \subseteq ... \subseteq H_m = \hat{H}$ which obeys the following conditions for all $i\in [m]$:

    \begin{enumerate}[label={\upshape (\arabic*)}]
        \item\label{case2:Hbasics} $A\cap V(H_i)\ne\emptyset$, $B\cap V(H_i)\ne\emptyset$,  and $H\cap H_i\ne\emptyset$, and $H_i$ is connected.
        \item\label{case2:pathsExist} Every component of $H_i - B$ contains a vertex of $A$,
        \item\label{case2:threemore} There are at most three more 
        edge-paths in 
        $Q\cup H_i$ than in $Q\cup H_{i-1}$,  
        \item\label{case2:onemore} 
        $Q\cup H_i$  has strictly more cross-inducing subgraphs than $Q\cup H_{i-1}$,
        \item\label{case2:eachx}  We set $H_i^*:=H_i-A$. For any vertex $x\in H_i^*$, there exists a connected subgraph $H_i^x$ of $H_i$ with $x\in H_i^x$ and $B\cap V(H_i^x)\ne\emptyset$ such that every component of $(H_i - H_i^x)- B=H_i-(H_i^x\cup B)$ contains a vertex of $A$,
        \item\label{case2:singletonB} If $x\in B\cap H_i^*$, then $H_i^x$ is the singleton $x$,
        \item\label{case2:eachxy} For any $x,y\in H_i^*$, if $x\ne y$ then $x\notin H_i^y$ or $y\notin H_i^x$.
    \end{enumerate}
    
    Using Lemma \ref{lm:bipartiteTruths} \ref{case:noncooperativeNeighbors}, there exists a $(\hat{H},X)$-path $P$ from some $b\in B$ to some $a\in A$ which contains some edge $e_1$ of $\hat{H}$. Let $H_1 = P$, and~\ref{case2:Hbasics} and ~\ref{case2:pathsExist} are clearly satisfied. The addition of $P$ to $Q$ can add at most three new edge-paths, so~\ref{case2:threemore} holds. 
    By Lemma \ref{lm:onecross}, $X\cup H_1$ has a cross-inducing subgraph which includes the path $P$. This cross-inducing subgraph contains the edge $e_1$ which does not belong to $Q$, so $Q\cup H_1$ has at least one more cross-inducing subgraph than $Q$. So \ref{case2:onemore} holds.
    For any $x\in H_1^*$, we define $H_1^x$ as the path $bPx$, which easily satisfies~\ref{case2:eachx} -~\ref{case2:eachxy}. 

    Suppose that for some $i\ge 1$ we have constructed $H_i$ as desired, and $H_i\ne \hat{H}$. 

    First we will show that there is a vertex $z_1\in H_i-N_H$ and edge $e_{i+1}\in E(\hat{H})-E(H_i)$ that has $z_1$ as one of its envertices. We will denote the other endvertex by $v$.

   Assume first that $H\neq H_i-N_H$
    Since $H$ is connected, this means $E(H)-E(H_i)$ is non-empty. As $H\cap H_i\ne\emptyset$, there must be some vertex $z_1\in H_i\cap H= H_i - N_H$ that is incident upon an edge $e_{i+1}\in E(H)-E(H_i)\subseteq E(\hat{H})-E(H_i)$. 
    
    Otherwise, $H=H_i-N_H$ and every edge in $E(\hat{H})-E(H_i)$ has one endpoint in $N_H$ and the other in $H$. In this case, we let $e_{i+1}$ be any edge in $E(\hat{H})-E(H_i)$ and $z_1$ its endpoint in $H$. 

    We achieved our first goal.

    Next we will find a $(\hat{H},H_i\cup X$)-path $R$ connecting 
    the vertex $z_1\in H_i-N_H=H_i^*-B$ to a vertex $z_2\in H_i\cup X$ that contains $e_{i+1}$. 
    
    If $v\in H_i\cup X$, simply let $R=e_{i+1}$ and $z_2=v$. If $v\notin H_i\cup X$,  by Lemma \ref{lm:bipartiteTruths} \ref{case:decomposition}
    we have $X\cup\hat{H}$ is $2$-connected, consequently there exists a pair of internally disjoint $(\hat{H}\cup X,H_i\cup X)$ paths from $v$ to distinct vertices of $H_i\cup X$. As at most one of these paths contain the vertex $z_1$, we let $P$ be a path connecting $v$ to $H_i\cup X$ which does not contain the vertex $z_1$ (and consequently does not contain the edge $e_{i+1}$). As in the proof of the previous lemma, $P$ is a $(\hat{H},H_i\cup X)$-path. Set $R=Pe_{i+1}$ and $z_2$ to be the endpoint of $P$ different from $v$.

    We achieved our second goal.

    We let $H_{i+1}=H_i\cup R$. Then~\ref{case2:Hbasics} holds immediately since $H_i\subseteq H_{i+1}$ and the endvertex $z_1$ of $R$ is in $H_i$. Every component of $H_{i+1}-B$ contains a component of $H_i-B$ unless both endpoints of $R$ lie in $B$, but this is impossible because $z_1\in H_i^*-B$. So~\ref{case2:pathsExist} holds. Furthermore, $R$ is an $H_{i+1}$-edge-path, so ~\ref{case2:threemore} holds. We must now show that $H_{i+1}$ satisfies conditions \ref{case2:onemore}-\ref{case2:eachxy}. Note that either $z_2\in A$, $z_2\in B$, or $z_2\in H_i^*-B$. 

    \textit{Case 1: $z_2\in A$:}
    
    As $H_i$ satisfies~\ref{case2:eachx}, $H_i^{z_1}$ contains a path $T$ from some $b\in B$ to $z_1$. Since $T\subseteq H_i^*$ and $R\cap H_i^*=\{z_1\}$, $bTz_1Rz_2$ is a $(\hat{H},X)$ path from $b$ to $z_2$ containing the edge $e_{i+1}$. By Lemma \ref{lm:onecross}, $(bTz_1Rz_2)\cup X$ contains a cross-inducing subgraph which contains $bTz_1Rz_2$ and thus contains the edge $e_{i+1}$. As $(bTz_1Rz_2)\cup X\subseteq H_{i+1}\cup Q$ and $H_i\cup Q$ does not contain $e_{i+1}$,
    \ref{case2:onemore} is satisfied. 
    
    Now, for every $x\in H_{i+1}^*$ we define $H_{i+1}^x$ and show these choices satisfy \ref{case2:eachx} - \ref{case2:eachxy}. 
    
    First, suppose $x\in H_i^*$. Then we set $H_{i+1}^x = H_i^x$.  This clearly preserves \ref{case2:singletonB}. Since $H_i^x$ satisfies \ref{case2:eachx}, $H_{i+1}^x$ is connected and contains $x$ and some vertex of $B$. Furthermore, any component of $(H_{i+1} - H_{i+1}^x)-B$  either contains $z_2$ or contains a component of $(H_{i} - H_{i}^x)-B$, so \ref{case2:eachx} is satisfied. 
    
   Next, if $x\in H_{i+1}^*-H_i^*=H_{i+1}-H_i$, then $x$ is an internal vertex of $R$. In particular, $x\notin B$, so we need not worry about \ref{case2:singletonB} in this instance. We set $H_{i+1}^x=xRz_1\cup H_i^{z_1}$. This is clearly connected and contains both $x$ and a vertex of $B$. Furthermore, any component of $(H_{i+1} - H_{i+1}^x)-B$  is either a component of $(H_{i} - H_{i}^{z_1})-B$  or contains $z_2$, so \ref{case2:eachx} is satisfied. 
   
   Now we check condition ~\ref{case2:eachxy}. Let $x, y\in H_{i+1}^*$, $x\ne y$. If both $x,y$ are internal vertices of $R$, then $H_{i+1}^x=H_i^x\cup xRz_1$ and $H_{i+1}^y=H_i^y\cup yRz_1$, so ~\ref{case2:eachxy} follows easily. Otherwise without loss of generality $x\in H_i^*$,. If $y\in H_i^*$, then ~\ref{case2:eachxy} follows from the induction hypothesis, otherwise $y\notin H_{i+1}^x$. Therefore ~\ref{case2:eachxy} holds, so we are done with the case where $z_2\in A$.

   \textit{Case 2: $z_2\in B$}

   As $H_i$ satisfies \ref{case2:pathsExist}, there exists a path $T$ from $z_1$ to some $a\in A$ with $T\subseteq H_i - B$. Since $V(R)\cap V(H_i-B) = \{z_1\}$, the path $z_2Rz_1Ta$ is a $(\hat{H},X)$ path from a vertex of $B$ to a vertex of $A$ containing the edge $e_{i+1}$. By Lemma \ref{lm:onecross}, $(z_2Rz_1Ta)\cup X$ contains a cross-inducing subgraph which contains $z_2Rz_1Ta$ and thus contains the edge $e_{i+1}$. As
   $(z_2Rz_1Ta)\cup X\subseteq H_{i+1}\cup Q$ and $H_i\cup Q$ does not contain $e_{i+1}$,
\ref{case2:onemore} is satisfied. 
   
   Now, for every $x\in H_{i+1}^*$ we define $H_{i+1}^x$ and show these choices satisfy \ref{case2:eachx} - \ref{case2:eachxy}.

    First we show \ref{case2:eachx} and \ref{case2:singletonB}:
    
    Suppose $x\in H_i^*$. We have four possibilities. 
    \begin{enumerate}[label={\upshape (\roman*)}]
        \item\label{subcase0} If $z_1, z_2\in H_i^x$, then we set $H_{i+1}^x = H_i^x \cup R$,
        \item\label{subcase1} If $z_1\in H_i^x$ and $z_2\notin H_i$, then we set $H_{i+1}^x = H_i^x \cup R$,
        \item\label{subcase2} If $z_1\in H_i^x$ and $z_2\in H_i - H_i^x$, then we set $H_{i+1}^x = H_i^x \cup (R- \{z_2\})$,
        \item\label{subcase3} If $z_1\notin H_i^x$,
        then we set $H_{i+1}^x = H_i^x$,
    \end{enumerate}

    In each case, it is clear that $H_{i+1}^x$ is connected and contains both $x$ and some vertex of $B$. Since $H_{i+1} = H_i\cup R$, it is not hard to see that in each case, every component of $(H_{i+1} - H_{i+1}^x)-B$  is either a component of $(H_i - H_i^x)-B$  or is the union of a component of $(H_i - H_i^x)-B$  that contains $z_1$ with $R$ (this happens only in case  \ref{subcase3}). So, \ref{case2:eachx} holds. 
    Assume $x\in B$, then $H_i^x=x$, and, since $z_1\notin B$, we must be in \ref{subcase3}\, therefore $H_{i}^x=\{x\}=H_{i+1}^x$. Thus, \ref{case2:singletonB} holds as well.

    Next, suppose $x\in H_{i+1}^* - H_i^*$. Then, $x\in R-z_1$.
    We set $H_{i+1}^x = xRz_2$. This is clearly connected, contains $x$, and $z_2\in B$;. Moreover, the only vertex of $B$ in $R$ is $z_2$, and $R$ contains no vertex of $A$. Therefore,
    if $x\in B$ then $x=z_2$ and $H_{i+1}^x$ is the path $z_2Rz_2$, which is the singleton $z_2$, so ~\ref{case2:singletonB} holds.
    As $(H_{i+1}-H_{i+1}^x)-B=((H_i\cup R)-(xRz_2))-B=(H_i-B)\cup (z_1Rx-x)$ and $z_1\in H_i-B$,
    the components of $(H_{i+1}-H_{i+1}^x)-B$ each contain a component of $H_i-B$. Since~\ref{case2:pathsExist} holds for $H_i$, \ref{case2:eachx} holds for $H_{i+1}$

    Finally we must show that \ref{case2:eachxy} holds with the choices we made. Let $x,y\in H_{i+1}^*$ be distinct. If $x,y\in H_{i}^*$, then condition \ref{case2:eachxy} is satisfied since the vertices of $H_{i+1}^x$ and $H_{i+1}^y$ differ from those of $H_{i}^x$ and $H_i^y$ only possibly by internal vertices of $R$. If $x,y$ are both internal vertices of $R$, then if $x\in H_{i+1}^y$ it must be that $x\in z_2Ry$. Therefore, $H_{i+1}^x = z_2Rx$ so $y\notin H_{i+1}^x$. So it suffices to assume that $x\in H_i^*$ and $y$ is an internal vertex of $R$. This means that $H_{i+1}^y = yRz_2$. If $x\in yRz_2\cap H_i^*$, then $x = z_2$.
    But since $z_2\in B$, by \ref{case2:singletonB} we have that $H_{i+1}^x = H_{i+1}^{z_2}$ is singleton $z_2$. So $y\notin H_{i+1}^x$. So \ref{case2:eachxy} holds and we are done.

    \textit{Case 3: $z_2\in H_i^*-B$:}
     
    As \ref{case2:eachx} and \ref{case2:eachxy} hold for $H_i$, we have subgraphs $H_i^{z_1}$ and $H_i^{z_2}$ of $H_i^*$ such that for some $j\in[2]$, $z_{3-j}\notin H_i^{z_{j}}$. By \ref{case2:eachx}, $H_i^{z_j}$ contains a path $T$ from $z_j$ to some $b\in B$, and every component of $(H_i - H_i^{z_j})-B$  contains some vertex of $A$. Let $S$ be a path in $(H_i - H_i^{z_j})-B$ from $z_{3-j}$ to some $a\in A$. Note that $V(S)\cap N_H = \{a\}$ and $V(T)\cap N_H = \{b\}$.   
    Now, set $P^* = bTz_jRz_{3-j}Sa$. This $P^*$ is a $(\hat{H},X)$ path from $b$ to $a$ containing the edge $e_{i+1}$. By Lemma \ref{lm:onecross}, $P^*\cup X$ contains a cross-inducing subgraph which contains $P^*$ and thus contains the edge $e_{i+1}$.
    As $P^*\cup X\subseteq Q\cup H_{i+1}$ and $Q\cup H_i$ does not contain the edge $e_{i+1}$, \ref{case2:onemore} is satisfied.

    Now, for every $x\in H_{i+1}^*$ we define $H_{i+1}^x$ and show these choices satisfy \ref{case2:eachx} - \ref{case2:eachxy}. 
    
    First, suppose $x\in H_i^*$. If $z_1,z_2\in H_i^x$, then we set $H_{i+1}^x = H_i^x \cup R$. Otherwise, we set $H_{i+1}^x = H_i^x$. In both cases $H_{i+1}^x$ is connected and contains $x$ and some vertex of $B$. As $H_{i+1}=H_i\cup R$, $H_{i+1}^x\in\{H_i^x,H_i^x\cup R\}$, and if $H_{i+1}^x=H_i^x$ then one of $z_1,z_2$ is not in $H_i^x$,
    every component of $(H_{i+1} - H_{i+1}^x)-B$ contains a component of $(H_{i} - H_{i}^x)-B$, and thus contains a vertex of $A$. So \ref{case2:eachx} holds in this case.
    Moreover, if $x\in B$, i.e. $x\in B\cap H_i^*$, then, as $H_i^x$ is the singleton $x$ and $z_2\notin B$, $H_{i+1}^x=H_i^x$ and 
    \ref{case2:singletonB} follows.

    Next, if $x\in H_{i+1}^* - H_i^*$, then $x$ is an internal vertex of the path $R$. In this case, $x\notin B$ so we need not consider \ref{case2:singletonB}. Recall that $j$ was chosen such that $z_{3-j}\notin H_i^{z_j}$. In this case, we set $H_{i+1}^x = H_i^{z_j} \cup xRz_j$. Clearly, $H_{i+1}^x$ is connected and contains $x$ and some vertex of $B$. Every component of $(H_{i+1} - H_{i+1}^x)-B$ is either a component of $(H_i - H_i^x)-B$, and thus contains a vertex in $A$, or contains the component of $(H_i - H_i^{z_j})-B$ containing $z_{3-j}$, and thus also contains some vertex of $A$.
    So \ref{case2:eachx} is satisfied for our choices of $H_{i+1}^x$.
    
    Finally we must show that \ref{case2:eachxy} holds. Let $x,y\in H_{i+1}^*$ be distinct. If $x,y\in H_{i}^*$, then condition \ref{case2:eachxy} is satisfied since the vertices of $H_{i+1}^x$ and $H_{i+1}^y$ differ from those of $H_{i}^x$ and $H_i^y$ only possibly by internal vertices of $R$. Furthermore, if $x,y$ are both internal vertices of $R$, then if $x\in H_{i+1}^y$ it must be that $x\in z_jRy$. Therefore, $y\notin H_{i+1}^x$. So it suffices to assume that $x\in H_i^*$ and $y$ is an internal vertex of $R$. Suppose $x\in H_{i+1}^y$. Then by construction, $x\in H_i^{z_j}$. Since $x\in H_i^*$, $H_{i+1}^x$ can only contain $y$ if $H_i^x$ contains both $z_1$ and $z_2$. But since $x\in H_i^{z_j}$, $z_j\notin H_i^x$. Therefore $y\notin H_{i+1}^x$, and we are done.

    Therefore, we have the desired sequence. Note that by~\ref{case2:onemore}, $Q\cup H_m = Q\cup \hat{H}$ contains at least $m$ more cross-inducing subgraphs than $Q$, and by~\ref{case2:threemore}, $Q\cup H_m$ has at most $3m$ more edge-paths than $Q$. Since $Q \cup H_m$ has $k'$ more cross inducing subgraphs than $Q$, $m \leq k'$. Therefore, $Q \cup H_m$ has at most $3m \leq 3k'$ more edge-paths than $Q$. This completes the proof.

\end{proof}

\begin{lemma}\label{lm:forkingLemma}
Let $w,a_1,a_2,a_3$ be four different vertices, for 
$i\in\{1,2,3\}$ let $T_i$ be a 
$w$-$a_i$ path such that $T_1,T_2,T_3$ 
are internally disjoint, and set $T=T_1\cup T_2\cup T_3$.

   Let $Y$ be a path between two vertices of $T$ such that $Y$ is internally disjoint from $T$ and $Y$ contains at most one of $a_1,a_2,a_3$. Then, for any  $x_1\in V(T)$ and $x_2\in V(Y)-V(T)$, there exists $u\in V(T\cup Y)-\{a_1,a_2,a_3\}$ such that 
   for each $i\in\{1,2,3\}$ the graph $T\cup Y$ contains a $u$-$a_i$ path $P_i$ where $P_1,P_2,P_3$ are internally disjoint, and at least one of the following holds:

    \begin{enumerate}[label={\upshape (\Alph*)}]
        \item\label{case:fork1} $x_1,x_2\in V(P_1\cup P_2\cup P_3)$,
        \item \label{case:fork2} $x_1,x_2\notin \{a_1,a_2,a_3\}$, and (potentially switching indices between $x_1$ and $x_2$) for 
        some $j\in[3]$ we have $x_1\in V(P_j)$, and there exists a $(T\cup Y, P_1\cup P_2\cup P_3)$ path $S$ with endpoints $s_1,s_2$ such that $x_2$ is an internal vertex of $S$ and one of the following is satisfied:
        \begin{enumerate}[label={\upshape (\Roman*)}]
            \item\label{case:fork3} $s_1\in V(uP_jx_1)-\{x_1\}$, $s_2\in V(x_1P_ja_j)-\{x_1\}$.
            \item\label{case:fork4}for some $k\in[3]-\{j\}$, $s_1\in V(P_k)-\{u,a_k\}$, $s_2 = a_j$.
        \end{enumerate}

    \end{enumerate}
\end{lemma}
\begin{figure}[t]
\begin{tikzpicture}[scale=0.8]

\begin{scope}[xshift=-10cm, yshift = 5cm]

\node at (-2,4) {\large $(i)$};

\node at (-2,0.75) { $T_1$};
\node at (-0.3,0.75) { $T_2$};
\node at (1.25,0.75) { $T_3$};

\node[vertexc] (a1) at (-2,0) {$a_1$};
\node[vertexc] (a2) at (0,0) {$a_2$};
\node[vertexc] (a3) at (2,0) {$a_3$};

\node[vertexc] (u) at (0,4) {$w$};

\draw[thickedge] (a1)--(u);
\draw[thickedge] (u)--(0,2);
\draw[thickedge] (0,2)--(a2);
\draw[thickedge] (a3)--(1.5,1);
\draw[thickedge] (0.5,3) -- (u);
\draw[thickedge, dashed] (0.5,3)--(1.5,1);

\node[vertex, color = black] (y2) at (0.5,3) {$.$};
\node[vertex, color = black] (y1) at (1.5,1) {$.$};

\node at (0.7,3.4) {\color{black} $y_1$};
\node at (1.9,0.9) {\color{black} $y_2$};

\node[vertex, color = gray] (z1) at (-1,2) {$.$};
\node at (-1.4,2) {\color{black} $x_1$};


\draw[rededge] (1.5,1)  to [out=40, in=0] (0.5,3);
\end{scope}

\begin{scope}[xshift=-5cm, yshift=5cm]
\node at (-2,4) {\large $(ii)$};
\node at (-2,0.75) { $T_1$};
\node at (-0.3,0.75) { $T_2$};
\node at (1.25,0.75) { $T_3$};

\node[vertexc] (a1) at (-2,0) {$a_1$};
\node[vertexc] (a2) at (0,0) {$a_2$};
\node[vertexc] (a3) at (2,0) {$a_3$};

\node[vertexc] (u) at (0,4) {$w$};

\draw[thickedge] (a1)--(u);
\draw[thickedge] (u)--(0,2);
\draw[thickedge] (0,2)--(a2);
\draw[thickedge] (a3)--(1,2);
\draw[thickedge, dashed] (u)--(1,2);

\node[vertex, color = black] (y2) at (1,2) {$.$};
\node[vertex, color = black] (y1) at (0,2) {$.$};

\node at (-0.35,2) {\color{black} $y_1$};
\node at (1.5,2) {\color{black} $y_2$};

\node[vertex, color = gray] (z1) at (-1,2) {$.$};
\node at (-1.4,2) {\color{black} $x_1$};


\draw[rededge] (y1) -- (y2);
\end{scope}

\begin{scope}[xshift=0cm, yshift = 5cm]
\node at (-2,4) {\large $(iii)$};
\node at (-2,0.75) { $T_1$};
\node at (-0.3,0.75) { $T_2$};
\node at (1.25,0.75) { $T_3$};

\node[vertexc, color=gray] (a1) at (-2,0) {\color{white}$a_1$};
\node[vertexc] (a2) at (0,0) {$a_2$};
\node[vertexc] (a3) at (2,0) {$a_3$};

\node[vertexc] (u) at (0,4) {$w$};

\draw[thickedge] (a3)--(u);
\draw[thickedge] (u)--(0,2);
\draw[thickedge] (0,2)--(a2);
\draw[thickedge] (a1)--(-1.5,1);
\draw[thickedge] (-0.5,3) -- (u);
\draw[thickedge, dashed] (-0.5,3)--(-1.5,1);

\node[vertex, color = black] (y2) at (-0.5,3) {$.$};
\node[vertex, color = black] (y1) at (-1.5,1) {$.$};

\node at (-0.7,3.4) {\color{black} $y_1$};
\node at (-1.1,0.9) {\color{black} $y_2$};

\node at (-1.4,0) {\color{black} $x_1$};


\draw[rededge] (y2)  to [out=-160, in=120] (y1);
\end{scope}

\begin{scope}[xshift=5cm, yshift=5cm]
\node at (-2,4) {\large $(iv)$};
\node at (-2,0.75) { $T_1$};
\node at (-0.3,0.75) { $T_2$};
\node at (1.25,0.75) { $T_3$};

\node[vertexc, color=gray] (a1) at (-2,0) {\color{white}$a_1$};
\node[vertexc] (a2) at (0,0) {$a_2$};
\node[vertexc] (a3) at (2,0) {$a_3$};

\node[vertexc] (u) at (0,4) {$w$};

\draw[thickedge] (a1)--(-1,2);
\draw[thickedge, dashed] (-1,2)--(u);
\draw[thickedge] (u)--(0,2);
\draw[thickedge] (0,2)--(a2);
\draw[thickedge] (a3)--(1,2);
\draw[thickedge] (u)--(1,2);

\node[vertex, color = black] (y1) at (-1,2) {$.$};
\node[vertex, color = black] (y2) at (0,2) {$.$};

\node at (0.35,2) {\color{black} $y_2$};
\node at (-1.5,2) {\color{black} $y_1$};

\node at (-1.4,0) {\color{black} $x_1$};


\draw[rededge] (y1) -- (y2);
\end{scope}

\begin{scope}[xshift=-10cm]
\node at (-2,3.5) {\large $(v)$};
\node at (-2,0.75) { $T_1$};
\node at (-0.3,0.75) { $T_2$};
\node at (1.25,0.75) { $T_3$};

\node[vertexc] (a1) at (-2,0) {$a_1$};
\node[vertexc] (a2) at (0,0) {$a_2$};
\node[vertexc] (a3) at (2,0) {$a_3$};

\node[vertexc] (u) at (0,4) {$w$};

\draw[thickedge] (a1)--(u);
\draw[thickedge, dashed] (u)--(0,1);
\draw[thickedge] (0,1)--(a2);
\draw[thickedge] (a3)--(u);

\node[vertex, color =black] (y1) at (0,1) {$.$};
\node[vertex, color = black] (y2) at (-.5,3) {$.$};

\node at (0.4,1) {\color{black} $y_2$};
\node at (-.9,3) {\color{black} $y_1$};

\node[vertex, color = gray] (z1) at (-1,2) {$.$};
\node at (-1.4,2) {\color{black} $x_1$};


\draw[rededge] (y1) -- (y2);
\end{scope}

\begin{scope}[xshift=-5cm]
\node at (-2,3.5) {\large $(vi)$};
\node at (-2,0.75) { $T_1$};
\node at (-0.3,0.75) { $T_2$};
\node at (1.25,0.75) { $T_3$};

\node[vertexc, color=black] (a1) at (-2,0) {\color{white}$a_1$};
\node[vertexc] (a2) at (0,0) {$a_2$};
\node[vertexc] (a3) at (2,0) {$a_3$};

\node[vertexc] (u) at (0,4) {$w$};

\draw[thickedge] (a1)--(u);
\draw[thickedge] (u)--(0,2);
\draw[thickedge] (0,2)--(a2);
\draw[thickedge] (a3)--(u);

\node[vertex, color = black] (y1) at (0,2) {$.$};

\node at (0.4,2) {\color{black} $y_2$};
\node at (-1.4,0) {\color{black} $y_1$};

\node[vertex, color = gray] (z1) at (-1,2) {$.$};
\node at (-1.4,2) {\color{black} $x_1$};


\draw[rededge] (0,2) -- (a1);
\end{scope}

\begin{scope}[xshift=0cm, yshift = 0cm]
\node at (-2,3.5) {\large $(vii)$};
\node at (-2,0.75) { $T_1$};
\node at (-0.3,0.75) { $T_2$};
\node at (1.25,0.75) { $T_3$};

\node[vertexc] (a1) at (-2,0) {$a_1$};
\node[vertexc] (a2) at (0,0) {$a_2$};
\node[vertexc] (a3) at (2,0) {$a_3$};

\node[vertexc] (u) at (0,4) {$w$};

\draw[thickedge] (a3)--(u);
\draw[thickedge] (u)--(0,2);
\draw[thickedge] (0,2)--(a2);
\draw[thickedge] (a1)--(-1,2);
\draw[thickedge] (-0.5,3) -- (u);
\draw[thickedge, dashed] (-0.5,3)--(-1,2);

\node[vertex, color = black] (y2) at (-0.5,3) {$.$};
\node[vertex, color = black] (y1) at (-1,2) {$.$};

\node at (-0.7,3.4) {$y_1$};
\node at (-0.6,2) {$y_2$};

\node[vertex, color = gray] (z1) at (-1.5,1) {$.$};
\node at (-1.1,1) {\color{black} $x_1$};


\draw[rededge] (y2)  to [out=-160, in=120] (y1);
\end{scope}

\begin{scope}[xshift=5cm]
\node at (-2,3.5) {\large $(viii)$};
\node at (-2,0.75) { $T_1$};
\node at (-0.3,0.75) { $T_2$};
\node at (1.25,0.75) { $T_3$};

\node[vertexc] (a1) at (-2,0) {$a_1$};
\node[vertexc] (a2) at (0,0) {$a_2$};
\node[vertexc] (a3) at (2,0) {$a_3$};

\node[vertexc] (u) at (0,4) {$w$};

\draw[thickedge] (a1)--(u);
\draw[thickedge] (u)--(0,2);
\draw[thickedge] (0,2)--(a2);
\draw[thickedge] (a3)--(u);

\node[vertex, color =black] (y2) at (-1.5,1) {$.$};
\node[vertex, color=black] (y1) at (-0.5,3) {$.$};

\node at (-0.7,3.4)  {\color{black} $y_1$};
\node at (-1.1,0.9) {\color{black} $y_2$};

\node[vertex, color = gray] (z1) at (-1,2) {$.$};
\node at (-1.4,2) {$x_1$};


\draw[rededge] (-1.5,1)  to [out=140, in=180] (-0.5,3);
\end{scope}
\end{tikzpicture} 
\caption{Selection of possible $T\cup Y$ in Lemma~\ref{lm:forkingLemma}. In each case the edge-paths of $T$ are drawn as straight lines from $w$ to the corresponding $a_j$. The path $Y$ is drawn as a thick curve/line between vertices of $T$, with labeled endpoints. In each case, $x_2$ belongs to this thick path, while a location of $x_1$ is shown explicitly as a gray vertex. The non-dashed parts of each drawing show the graph $P_1\cup P_2\cup P_3$ in case \ref{case:fork1} and $P_1\cup P_2\cup P_3\cup S$ in case \ref{case:fork2}. We note that case \ref{case:fork2} is represented only in drawings (vi) and (viii), and all other cases represent case \ref{case:fork1}.}
\label{fig:casesForTripleForks}
\end{figure}
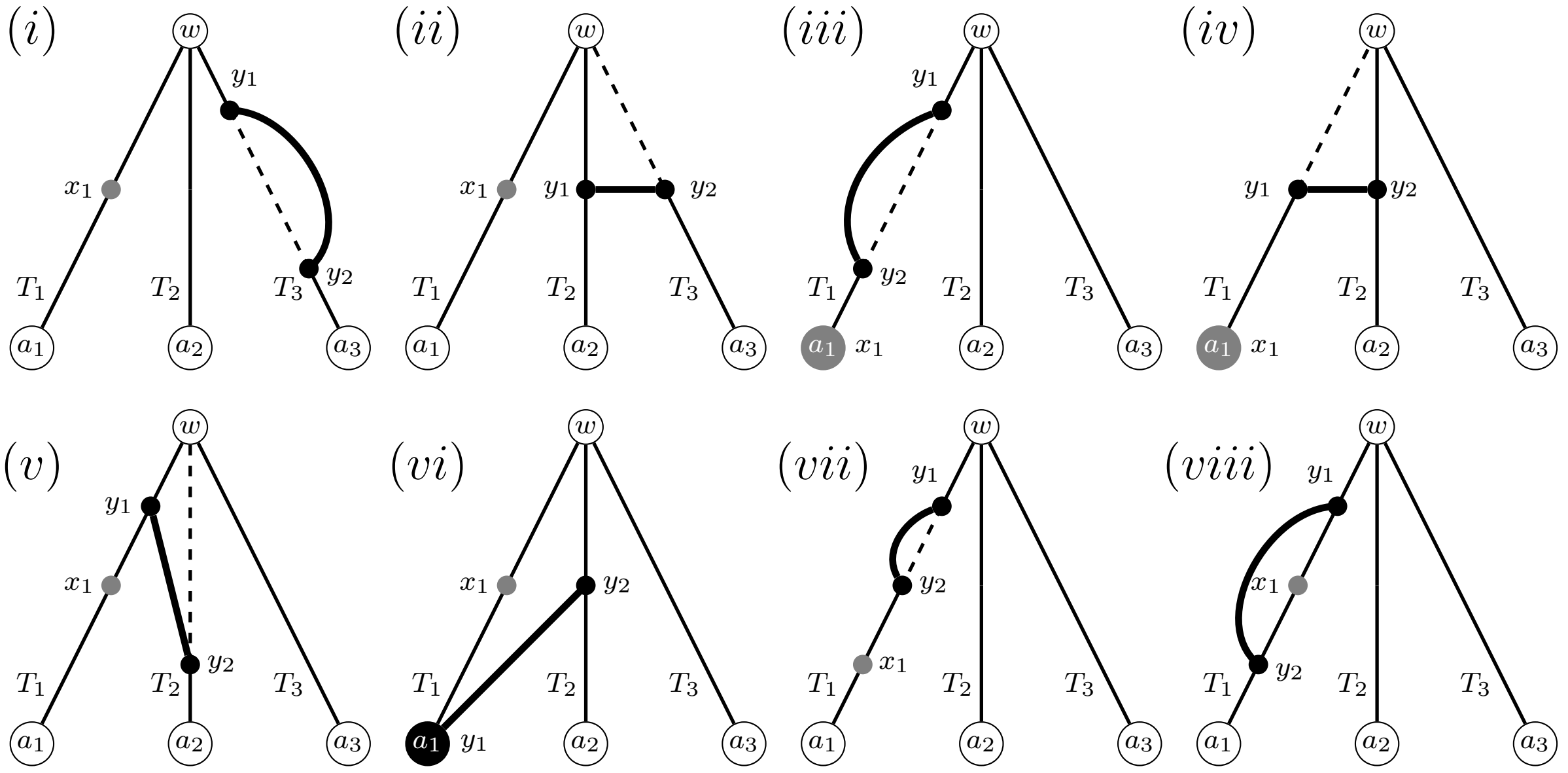

\begin{proof}
   $w,a_1,a_2,a_3,T_1,T_2,T_3,T$ and $Y$ are given. Let $y_1,y_2$ denote the endpoints of $Y$, and let $x_1\in V(T)$ and $x_2\in V(Y)-V(T)$ be arbitrary. Notice that $x_2\notin\{w,a_1,a_2,a_3\}$. Without loss of generality (we may switch the indeces of $a_1,a_2,a_3$) we may assume that $x_{1}\in V(T_1)$. We must find paths $P_1,P_2,P_3$ to satisfy the conclusion of the lemma.

   {\it Case 1} $y_1,y_2\in T_2\cup T_3$.
   
   If $y_1,y_2\in T_k$ for some $k\in\{2,3\}$, set $P_1=T_1$, $P_{5-k}=T_{5-k}$, and obtain $P_k$ from $T_k$ by replacing the sub-path  $y_1T_ky_2$ with $Y$  (see Figure~\ref{fig:casesForTripleForks} $(i)$). This clearly satisfies \ref{case:fork1}. Otherwise without loss of generality $y_1\in V(T_2)\setminus\{w\}$ and $y_2\in V(T_3)\setminus \{w\}$ with $y_1\ne a_2$. Set $P_1=y_1T_2wT_1a_1$,  $P_2=y_1T_2a_2$, $P_3=y_1Yy_2T_3a_3$  (see Figure~\ref{fig:casesForTripleForks} $(ii)$). These paths satisfy \ref{case:fork1}.

   {\it Case 2:} $\{y_1,y_2\}\not\subseteq V(T_2)\cup V(T_3)$. 
   
   Without loss of generality, $y_1\in V(T_1)\setminus\{w\}$ and $y_2\in V(T_1\cup T_2)$. In particular, either $y_2\in V(T_1)$ or $y_2\in V(T_2)\setminus\{w\}$. 

    {\it Subcase 2a:} $x_1=a_1$. 
    
    If $y_2\in V(T_1)$,  set $P_2=T_2, P_3=T_3$, and obtain $P_1$ from $T_1$ by replacing $y_1T_1y_2$ with $Y$ (see Figure~\ref{fig:casesForTripleForks}  $(iii)$) . If $y_2\in V(T_2)-\{w\}$, then we set $P_1 = y_2Yy_1T_1a_1, P_2 = y_2T_2a_2, P_3 = y_2T_2wT_3a_3$ (see Figure~\ref{fig:casesForTripleForks} $(iv)$).  These paths satisfy \ref{case:fork1}. 
    
    {\it Subcase 2b}: $x_1\ne a_1$. In particular $x_1,x_2\notin\{a_1,a_2,a_3\}$.

    Assume first that $y_2\in V(T_2)-\{w\}$. If $y_1\ne a_1$, set $P_1 = y_1T_1a_1$, $P_2=y_1Yy_2T_2a_2$, and $P_3=y_1T_1wT_3a_3$ (see Figure~\ref{fig:casesForTripleForks}$(v)$). Here, $u=y_1$, and these paths satisfy \ref{case:fork1}. If $y_1 = a_1$, then by assumption $y_2\ne a_2$. We set each $P_j=T_j$, and $S=Y$, and we satisfy \ref{case:fork2}\ref{case:fork4} with $j=1,k=2$ (see Figure~\ref{fig:casesForTripleForks} $(vi)$).
    
    Assume next that $y_2\in V(T_1)$, so $y_1,y_2\in V(T_1)$.
    If $y_1,y_2\in wT_1x_1$ or $y_1,y_2\in x_1T_1a_1$, then we set $P_2 = T_2, P_3=T_3$, and $P_1$ is formed by replacing the $y_1T_1y_2$ sub-path  in $T_1$ with $Y$ (see Figure~\ref{fig:casesForTripleForks} $(vii)$). These choices satisfy \ref{case:fork1} Otherwise we may assume that $y_1\in wT_1x_1, y_2\in x_1T_1a_j$, and $x_1\notin\{y_1,y_2\}$. But then, we set each $P_j=T_j$ and $S=Y$ and we satisfy \ref{case:fork2}\ref{case:fork3} (see Figure~\ref{fig:casesForTripleForks} $(viii)$).   
\end{proof}

\begin{lemma}\label{lm:k33LinearWithCooperative}
    Suppose $G$ is a 2-connected graph, containing no $K_5$ subdivisions, $X$ be a cross-inducing subgraph of $G$, and $H\in\mathcal{P}_2^X$ such that $N_H$ is exactly one of the $X$-main vertex classes. Then for any subgraph $Q$ of $G$ with $X\subseteq Q$ and $V(Q)\cap V(H) =\emptyset$, if $Q\cup \hat{H}$ 
    has exactly $k'$ more cross-inducing subgraphs than $Q$, then $Q\cup \hat{H}$ has at most $3k'$ more edge-paths than $Q$.
\end{lemma}
\begin{proof}

    Let $G,X,Q$ be as described, and say that $N_H = \{a_1,a_2,a_3\}$. We proceed by constructing a sequence of graphs $\emptyset = H_0 \subseteq H_1 \subseteq ... \subseteq H_m = \hat{H}$ obeying the following conditions for all $i\in [m]$.

    \begin{enumerate}[label={\upshape (\arabic*)}]
        \item \label{case3:connectedBasic} $H_i\cap H$ is connected,
        \item \label{case3:connected} $X\cup H_i$ is 2-connected,
        \item \label{case3:morePaths} There are at most three more $Q\cup H_{i}$ edge-paths than there are $Q\cup H_{i-1}$ edge-paths,
         \item \label{case3:moreInducing} $Q\cup H_{i}$ has strictly more cross-inducing subgraphs than $Q\cup H_{i-1}$,
        \item \label{case3:tripleConditionFork} For every $x_1,x_2\in V(H_i)$, there exists a $u\in V(H_i)-N_H$ and three internally disjoint paths $P_1,P_2,P_3$ in $H_i$ such that each $P_j$ has endpoints $u$ and $a_j$. At least one of the following holds:
        \begin{enumerate}[label={\upshape (\Alph*)}]
            \item \label{case3:tripleConditionFork1} $x_{1},x_2\in V(P_1\cup P_2\cup P_3)$.          
            \item \label{case3:tripleConditionFork2} $x_1,x_2\notin \{a_1,a_2,a_3\}$, and (potentially switching indices of $x_1$ and $x_2$) $x_1\in V(P_j)$ for some $j\in[3]$, and there exists a $(H_i, P_1\cup P_2\cup P_3)$ path $S$ with endpoints $s_1,s_2$ such that $x_2$ is an internal vertex of $S$ and one of the following is satisfied:
        \begin{enumerate}[label={\upshape (\Roman*)}]
            \item\label{case3:tripleConditionFork3} $s_1\in V(uP_jx_1)-\{x_1\}$, $s_2\in V(x_1P_ja_j)-\{x_1\}$,
            \item\label{case3:tripleConditionFork4} for some $k\in[3]-\{j\}$, $s_1\in V(P_k)-\{u,a_k\}$, $s_2 = a_j$.
        \end{enumerate}
  \end{enumerate}
    \end{enumerate}

    Using Lemma \ref{lm:bipartiteTruths} \ref{case:mainClassNeighbors}, 
    there are internally disjoint paths $P_1,P_2,P_3$ from some $u\in \hat{H}- N_H$ to $a_1,a_2,a_3$ respectively;     we set $H_1=P_1\cup P_2\cup P_3$. This clearly satisfies all the criteria.

    Suppose we have constructed $H_i$ as desired and $H_i\ne\hat{H}$. We will find an $(\hat{H}, H_i)$ path $R$ with endpoints $z_1\in H_i- N_H$, $z_2\in H_i$. By connectivity of $\hat{H}$, and the fact that $\hat{H}$ contains no edges with both endpoints in $N_H$, there is some edge $e\in E(\hat{H}) - E(H_i)$ with one endpoint $z_1\in H_i- N_H$ with other endpoint $v$. 
    If $v\in H_i$, then set $R=e$ and set $z_2 = v$. Otherwise, $v\notin H_i$. By Lemma \ref{lm:bipartiteTruths} \ref{case:decomposition}, we have $X\cup\hat{H}$ is $2$-connected, consequently there exists a pair of internally disjoint $(X\cup \hat{H},H_i)$ paths from $v$ to distinct vertices of $H_i$; let $P$ be one of these two paths that does not contain $z_1$ (and consequently does not contain $e$). Since for $i\ge 1$ we have $N_H\subseteq H_i$, $N_H$ separates $X$ and $ \hat{H}$,
    and $z_1\in H_i- N_H$, $P$ is a
    $(\hat{H},H_i)$-path.
    We denote $R = Pe$, and call by $z_2$ the endpoint of $R$ other than $z_1$.

    We set $H_{i+1} = H_i \cup R$. Note that $H_{i+1}\cap H = (H_i\cap H)\cup (R\cap H)$. By assumption, $H_i\cap H$ is connected. Since $R$ is a path containing an edge with at most one vertex outside of $H$, $R\cap H$ is connected. Since $z_1\in H\cap H_i\cap R$, $H_{i+1}\cap H$ is connected. So \ref{case3:connectedBasic} holds.
    Conditions \ref{case3:connected} and \ref{case3:morePaths} clearly hold for $H_{i+1}$. Since $H_i$ satisfies \ref{case3:tripleConditionFork} and $z_1\in H_i- N_H$, $z_2\in H_{i}$, for $\{x_1,x_2\}=\{z_1,z_2\}$ there exists a $u\in V(H_i- N_H)$ and $P_1,P_2,P_3$ internally disjoint paths with $P_j$ a $u-a_j$ path, satisfying \ref{case3:tripleConditionFork1}, or \ref{case3:tripleConditionFork2}\ref{case3:tripleConditionFork3}, or \ref{case3:tripleConditionFork2}\ref{case3:tripleConditionFork4}.  
   In every cases, as we may freely re-index the $a_j$ and $P_j$, without loss of generality we may assume that either $z_1\in V(P_1)$ or $z_1\notin V(P_1\cup P_2\cup P_3)$. Furthermore, in the \ref{case3:tripleConditionFork2} cases, we may re-index the $z_i$ such that $z_1\in V(P_1)$ since in this case both $z_i\in H_i\setminus N_H$. 

    First, we prove that in each case we can define pairwise internally disjoint paths $T_1,T_2,T_3$ where for some $w\in V(H_{i+1})-N_H$, each $T_j$ is a $w-a_j$ path, and $R\subseteq T_1\cup T_2\cup T_3$.

    \textit{Case: \ref{case3:tripleConditionFork1} with $z_2\in V(P_1)$}
    
    In this case, we define $T_1$ to be the path formed from $P_1$ by replacing $z_1P_1z_2$ with $R$. This is clearly a $u-a_1$ path. Then, define $T_2=P_2, T_3=P_3$. The $T_j$ paths are pairwise internally disjoint $u-a_j$ paths. $R\subseteq T_1$, and here $w=u$.

    \textit{Case: \ref{case3:tripleConditionFork1} with $z_2\notin V(P_1)$}
    
    Without loss of generality, we may assume that $z_2\in V(P_2)$.   
    In this case, we define the following paths: $T_1 = z_1P_1a_1$, $T_2=z_1Rz_2P_2a_2$, $T_3 = z_1P_1uP_3a_3$. Note that if $z_2=a_2$, then $T_2 = R$. The $T_j$ paths are pairwise internally disjoint $z_1-a_j$ paths. $R\subseteq T_2$, and here $w=z_1$.

    \textit{Case: \ref{case3:tripleConditionFork2}\ref{case3:tripleConditionFork3}}
    
    In this case, we define the following paths: $T_1 = uP_1z_1Rz_2Ss_2P_1a_1$, $T_2=P_2$, $T_3 = P_3$. The $T_j$ paths are pairwise internally disjoint $u-a_j$ paths. $R\subseteq T_1$, and here $w=u$.

    \textit{Case: \ref{case3:tripleConditionFork2}\ref{case3:tripleConditionFork4}}
    
    In this case, without loss of generality, $k=2$. We define the following paths: $T_1 = z_1P_1a_1$, $T_2=z_1Rz_2Ss_1P_2a_2$, $T_3 = z_1P_1uP_3a_3$. The $T_j$ paths are pairwise internally disjoint $z_1-a_j$ paths. $R\subseteq T_2$, and here $w=z_1$.

    Now, in every case, we have defined pairwise internally disjoint paths $T_1,T_2,T_3$ where for some $w\in V(H_{i+1})-N_H$, each $T_j$ is a $w-a_j$ path, and $R\subseteq T_1\cup T_2\cup T_3$. Set $T = T_1\cup T_2\cup T_3$. By Lemma \ref{lm:k331}, we have a cross-inducing subgraph of $Q\cup H_{i+1}$ containing $T$, and in particular containing $e$. So \ref{case3:moreInducing} holds.

    Now we must show \ref{case3:tripleConditionFork}. Suppose $x_1,x_2\in H_{i+1}$. If $x_1,x_2\in V(H_i)$ we have \ref{case3:tripleConditionFork} by assumption. If $x_1,x_2\in V(T)$, then $T_1,T_2,T_3$ satisfies \ref{case3:tripleConditionFork} \ref{case3:tripleConditionFork1}. So we may assume $x_1\in V(H_{i+1}-H_i)\subseteq V(R)\subseteq V(T)$, and $x_2\in V(H_{i+1}-T)$. 
    In particular, $x_1,x_2\notin \{a_1,a_2,a_3\}$.
    
    We will now find a $(H_{i+1}, T)$ path $Y$ with endpoints in $T$ that contains $x_2$ such that $Y$ contains at most one vertex of $N_H$. Once we do, then we are done by Lemma \ref{lm:forkingLemma}. Since $H_{i+1}\cup X$ is $2$-connected, $x_2\in V(H_{i+1}-T)$, $T\subseteq H_{i+1}$, and $N_H$ is a separating set for $H_{i+1}$ and $X$, there exists a $(H_{i+1}, T)$ path $Y'$ with endpoints in $T$ and $x_2\in V(Y')$. If $Y'$ has at most one endpoint in $N_H$, we use $Y=Y'$ and we are done. Otherwise, since $H\cap H_i$ is connected,  $T-N_H\subseteq H$, $z_1\in V(H)$, and $V(Y'-N_H)\cap V(T)=\emptyset$, there exists a $(H_{i+1}\cap H,T\cup Y')$ path $Y''$ with one endpoint in $Y'-N_H$ and other endpoint in $T-N_H$. $Y'\cup Y''$ contains the desired path $Y$.

\end{proof}

With these lemmas, we now proceed to the proof Lemma \ref{lm:k33Linear}.

\begin{proof}
     We enumerate the elements of $\mathcal{P}_2^X$ as $H_1,...,H_m$, and the edges of $F_X^*$ as $f_1,...,f_n$. We construct a sequence of graphs $X = G_0 \subseteq G_1 \subseteq ... \subseteq G_{m+n} = G(X)$ where for $i\in [n]$, $G_i = G_{i-1}\cup f_i$ and for $j\in [n+m]-[n]$, $G_j = G_{j-1}\cup \hat{H}_{j-n}$. For every $t\in [n+m]$, set $k_t$ to be the difference in the number of cross-inducing subgraphs of $G_{t}$ and $G_{t-1}$, and $p_t$ the difference in the number of edge-paths of $G_{t}$ and $G_{t-1}$. Clearly, $1 + \sum\limits_{t\in[m+n]} k_t \le k$ and $G(X)$ has $9+\sum\limits_{t\in[m+n]} p_t$ edge-paths.

    By Lemma \ref{lm:edgeestimate0} \ref{case:outsideEdges}, for every $i\in [n]$ we have that $G_i$ has at least one more cross-inducing subgraph than $G_{i-1}$ does, and furthermore $G_i$ has at most three more edge-paths than does $G_{i-1}$. In particular, $p_i\leq 3\le 3 k_i$ for $i\in [n]$. For $j\in [n+m]-[n]$, we have that $X\subseteq G_{j-1}$ and that $H_{j-n}\cap G_{j-1} = \emptyset$. By definition, either $N_{H_{j-n}}$ contains an $X$-subdividing vertex, contains two $X$-main vertices that are not $X$-cooperative, or is exactly an $X$-main vertex class. Therefore either Lemma \ref{lm:k33LinearwithSubdivision}, Lemma \ref{lm:k33LinearwithNoncooperative}, or Lemma \ref{lm:k33LinearWithCooperative} applies using $Q = G_{j-1}$. We have that $G_j$ has at most $3k_j$ more edge-paths than $G_{j-1}$. In particular, $p_j\leq 3k_j$ for $j\in [n+m]-[n]$.

    So the number of edge-paths of $G(X)$ is $$9 + \sum\limits_{t\in[m+n]} p_t \leq 9 + \sum\limits_{t\in[m+n]} 3k_t \le 9 + 3(k-1) = 6 + 3k$$
\end{proof}

\section{Proof of the main Theorems}\label{sec:mainproof}

Since $f(k)\ge k$, to show that for $k\le 2$ we have  $f(k)=k$, it is  enough to show that $f(k)\le k$

We will show first that $f(k)\le k$ for $k\le 1$ and $f(k)\le\binom{3k+7}{2}$ when $k\ge 2$.

We will use the following fact
\begin{claim}\label{cl:compute}
Let $s,\ell$ be nonnegative integers such that $s+\ell>0$,  and for each $i:1\le i\le \ell$ let $k_i\ge2$ be an integer. Set $m=s+\sum_{i=1}^{\ell} k_i$ Then
$$s+\sum_{i=1}^{\ell}\binom{3k_i+7}{2}\le\binom{3m+7}{2}
$$
Consequently, the function $h:\mathbb{N}\rightarrow\mathbb{N}$ defined as $h(k)=k$ when $k\le 1$ and $h(k)=\binom{3k+7}{2}$ when $k\ge 2$ is superadditive.
\end{claim}

\begin{proof}
Set $q=\sum_i k_i$. Then $m=s+q$.

Assume first that $\ell=0$. Then $1\le s$,
so $3s+6\ge 9$,
and $\binom{3s+7}{2}\geq\frac{9}{2}\cdot(3s+7)>s$. Therefore
the claim is true when $\ell=0$

Assume next that $\ell>0$. Then
\begin{eqnarray*}
\sum_{i=1}^{\ell}\binom{3k_i+7}{2}
&=&\frac{\sum_{i=1}^{\ell}\left(9k_i^2+39k_i+42\right)}{2}
\\
&=&\frac{9q^2-9\left(\sum_{i=1}^{\ell}\sum_{j\ne i}k_ik_j\right)+39q+42+42(\ell-1)}{2}\\
&=&\binom{3q+7}{2}-
\frac{9\left(\sum_{i=1}^{\ell}\sum_{j\ne i}k_ik_j\right)-42(\ell-1)}{2}\\
&\le&\binom{3q+7}{2}-\frac{36\ell(\ell-1)-42(\ell-1)}{2}\\
&=&\binom{3q+7}{2}-(\ell-1)(18\ell-21)
\le\binom{3q+7}{2},
\end{eqnarray*}
where the last inequality holds because either $\ell-1=0$ or $\ell-1\ge 1$ and $18\ell-21\ge 36-21>0$.
 In particular, the statement is true when $s=0$.

Now, we have that for $s,q>0$
\begin{eqnarray*}
\binom{3(s+q)+7}{2}
-\binom{3q+7}{2}
&=&\frac{(3s+3q+7)(3s+3q+6)-(3q+7)(3q+6)}{2}\\
&=&\frac{3s(3s+6q+17)}{2}>\frac{9s^2}{2}> s
\end{eqnarray*}
 
Using this, when $s,\ell>0$, we get
\begin{eqnarray*}
s+\sum_{i=1}^{\ell}\binom{3k_i+7}{2}
&\le&s+\binom{3q+7}{2}<
\binom{3(s+q)+7}{2}=\binom{3m+7}{2}
\end{eqnarray*}.
\end{proof}

For the rest of the proof, we consider a $2$-connected graph $G$ and suppress degree $2$ vertices as subdivided vertices (as we can insert them back in the drawing without difficulty). Thus, without loss of generality the edge-paths of $G$ are the edges of $G$ and every vertex of $G$ is a main vertex. Let $k$ be the number of cross-inducing subgraphs of $G$.

Assume first that $G$ is a $2$-connected graph that contains a $K_5$, then $k\ge 1$.
By Lemma~\ref{lm:K5linear}, $\X_G=G$ and $G$ has at most $3k+7$ edges.

If $k=1$, $G=\X_G$ is a $K_5$, 
consequently $G$ has a drawing with exactly one crossing. If $k>1$, put the main vertices of $G$ on a circle, connect any two vertices that has an edge connecting them  by a straight edge, and draw any additional edges sufficiently close to the first copy so no edges intersect more than once. As the crossing number of this drawing is at most the number of unordered pairs of edges, which is at most
$\binom{3k+7}{2}$, if $G$ has a subdivided $K_5$, it has a drawing with at most $h(k)$
crossings.

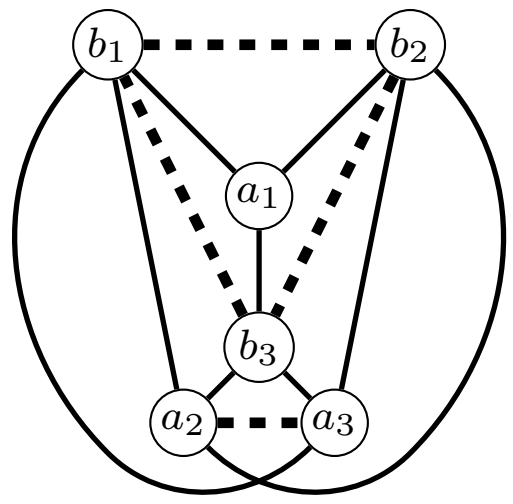
\begin{figure}[http]
\center
\begin{tikzpicture}
\node[vertexc] (b1) at (-1,2) {$b_1$};
\node[vertexc] (b2) at (1,2) {$b_2$};
\node[vertexc] (b3) at (0,0) {$b_3$};
\node[vertexc] (a1) at (0,1) {$a_1$};
\node[vertexc] (a2) at (-.5,-.5) {$a_2$};
\node[vertexc] (a3) at (.5,-.5) {$a_3$};
\draw[blueedge] (b1)--(a1)--(b2);
\draw[blueedge] (b2) to [out=-45,in=45] (1,-.7) to [out=225, in=-45] (a2);
\draw[blueedge] (a2)--(b3)--(a3);
\draw[blueedge] (b1) to [out=225,in=135] (-1,-.7) to [out=-45, in=225] (a3);
\draw[blueedge] (b1)--(a2);
\draw[blueedge] (b2)--(a3);
\draw[blueedge] (b3)--(a1);
\draw[rededge, dashed] (a2)--(a3);
\draw[rededge,dashed] (b1)--(b2)--(b3)--(b1);
\end{tikzpicture} 
\caption{A drawing of $G(X) \cup S^*$ when $k=1$. Edges of $S^*$ are thick and dashed.}
\label{fig:alternates}
\end{figure}

Now assume that $G$ contains no subdivided $K_5$, and there is a cross-inducing subgraph $X$ of $G$
such that $G=G(X)\cup S$, where $S$ is a 
set of edges connecting $X$-cooperative vertices of $X$ without parallel edges. 
Let $A=\{a_1,a_2,a_3\}$ and $B=\{b_1,b_2,b_3\}$ be the main vertex classes of $X$. 
Since $G=G(X)\cup S$ contains no subdivided $K_5$, by Lemma~\ref{lm:k331}~\ref{case:k5generate}
if $S$ contains more than $2$ edges, then $S$ is the set of edges connecting the $3$ vertices within the same main vertex class.
Therefore without loss of generality $S\subseteq S^*$ where $S^*$ contains an $a_2$-$a_3$ edge and the $3$ edges connecting the vertices in $B$. Let $k'$ be the number of cross-inducing subgraphs of $G(X)$.

As $\X_{G(X)}=G(X)\subseteq \X_G$, $1\le k'\le k$.

If $k=1$, then $k'=1$, and $G(X)=X$. 
Figure~\ref{fig:alternates} gives a drawing of $G(X)\cup S^*$ with $1=h(1)$ crossing such that only edges of $G(X)$ cross; this contains a drawing of $G(X)\cup S$ where only edges of $G(X)$ cross.

\begin{figure}[http]
\center
\begin{tikzpicture}
\node[vertexb] (a1) at (0,0) {};
\node[vertexb] (a2) at (1,0) {};
\node[vertexb] (a3) at (2,0) {};
\node[vertexb] (b1) at (0,-1) {};
\node[vertexb] (b2) at (1,-1) {};
\node[vertexb] (b3) at (2,-1) {};
\node[vertexb] (x) at (.2,-.5) {};
\node[vertexb] (y) at (.5,-.9) {};
\node[vertexb] (z) at (1.115,-.25) {};
\draw[blueedge] (a1)--(x)--(b1)--(a2)--(z)--(b2)--(a3)--(b3)--(a1)--(b2);
\draw[blueedge] (a2)--(b3);
\draw[blueedge] (a3)--(b1);
\draw[blueedge] (b2)--(y)--(x)--(z)--(y);
\draw[rededge,dashed] (a1)--(a2)--(a3) to [out=90,in=90] (a1);
\draw[rededge, dashed] (b1)--(b2)--(b3) to [out=-90,in=-90] (b1);
\end{tikzpicture}
\caption{Schematic sketch of the drawing of $G(X)\cup S$. Edges of $S$ are thick and dashed.}
\label{fig:gstar}
\end{figure}
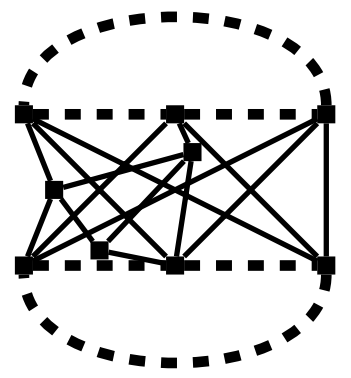

If $k\ge 2$, we create the following drawing:
take a rectangle, and place the vertices in the same main vertex class of $X$ on one side of the rectangle so the two main vertex classes are on two parallel lines. Place the remaining vertices of $G(X)$ inside the rectangle such that if three vertices of $G(X)$ are collinear, they form a main vertex class of $X$. Draw a straight line between two vertices of $G(X)$ if there is an edge between them in $G(X)$, add the necessary duplicate edges drawn close enough to their first copy so that only edges of $G(X)$ cross, and they cross only once. By Lemma~\ref{lm:k33Linear}, this
is a drawing of $G(X)$ with at most $\binom{3k'+7}{2}=h(k')$ crossings.
Now for each pair of $X$-cooperative vertices, draw
an edge connecting the two vertices either as a straight line or as a curved edge so the new edges do not cross each other and do not enter the rectangle
 (see Figure \ref{fig:gstar}). 
 This contains a drawing of $G(X)\cup S$ with at most $h(k')$ crossings such that only edges of $G(X)$ cross.

 Now Claim~\ref{cl:compute} and Lemma~\ref{lm:meta} gives that $f(k)\le h(k)$, which proves Theorem~\ref{th:existence}
and shows that Theorem~\ref{th:geza} is true for $k\ne 2$.

Now to finish the proof of Theorem~\ref{th:geza}, we need to show that
$f(2)=2$. Lemma~\ref{lm:meta}  gives that to prove $f(2)=2$, it is enough to create a drawing of the $2$-connected graphs $G$ with at most $2$ cross-inducing subgraphs in the cases when $G$ either contains a subdivided $K_5$, or $G$ does not contain a subdivided $K_5$ and
$G=G(X)\cup S$ for some cross-inducing subgraph $X$. As the
case when such graphs contain exactly $1$ cross-inducing subgraphs was done earlier, we may also assume that $G$ contains exactly $2$ cross-inducing subgraphs.

Assume first that $G$ contains a subdivided $K_5$, call it $X$. By Lemma~\ref{lm:K5linear} we have $\X_G=G$.
As $G$ contains another
cross inducing subgraph, $G$ has an edge $e$ not in $X$. As $G$ is $2$-connected, there is a $(G,X)$-path $P$ containing $e$; call the endvertices of $P$ $p$ and $q$.
By Lemma~\ref{lm:onecross}, if both $p,q$ are subdividing vertices of $X$, then $p,q$ are on the same edge-path of $X$, and $X\cup P$ has a cross-inducing subgraph containing $P$. As $X\cup P$ has $2$ cross-inducing subgraphs, and $G=\X_G=X\cup P$.
If $p,q$ are on the same edge-path of $X$,
take a drawing of $X$ with one crossing such that the edge-path containing $p,q$ is not crossed, and draw the path $P$ close to
the edge-path of $P$ so that it does not cross any edges.
Otherwise, the drawing on the left in Figure~\ref{fig:geza1} provides a drawing of $G$ with one crossing.

Now assume that $G$ contains no subdivided $K_5$, $X$ is a cross inducing subgraph of 
$G$ with $\X_G=X(G)$, and $G=G(X)\cup S$ where
$S$ is a set (not multiset!) of edges connecting $X$-cooperative vertices. We denote the main vertex classes of $G$ by
$A=\{a_1,a_2,a_3\}$ and $B=\{b_1,b_2,b_3\}$.

Since $X\subseteq G(X)\cup S$,
$G(X)$ contains at least $1$ cross-inducing subgraph.
If $G(X)$ contains precisely one cross-inducing subgraph, then $X=G(X)$, and, as $G(X)\cup S$ contains no $K_5$, the graph on Figure~\ref{fig:alternates} contains $G=G(X)\cup S$ has the requires drawing with fewer than $2$ crossings. 

Thus, we may assume that $G(X)$ contains two cross-inducing subgraphs, consequently
$X\subsetneq G(X)\cup$. We will first show that 
there is a $(G,X)$-path $P$ such that
$G(X)=X\cup P$.

As $G(X)$ contains exactly two cross-inducing subgraphs, Lemma~\ref{lm:edgeestimate0}~\ref{case:outsideEdges}
together with Lemma~\ref{lm:bipartiteTruths}~\ref{case:decomposition} gives that $G(X)=X\cup \{e\}$ where $e\in F^*_X$
or $G(X)=X\cup\hat{H}$ where $H\in\mathcal{P}_2^X$.
In the first case $P=e$.
If $G(X)=X\cup\hat{H}$ for some $H\in\mathcal{P}_2^X$,
Lemma~\ref{lm:bipartiteTruths} gives that $|N_H|=2$, and
by definition of $\mathcal{P}_2^X$ we have that the two vertices of $N_H$ are not $X$-cooperative.
As $\hat{H}$ is connected, it has a path $P$ between
the vertices of $N_H$. By Lemma~\ref{lm:bipartiteTruths} and the fact that $X\cup P$ contains at least two cross-inducing subgraphs, we have that $\hat{H}=P$

Now, we have that $G(X)=X\cup P$ where $P$ is a $(G,X)$-path whose endvertices are not $X$-cooperative and $\X_G=G(X)=X\cup P$ contains exactly two cross-inducing subgraphs. Lemma~\ref{lm:onecross} gives that either the endvertices of $P$ lie on the same edge-path of $X$, or one of them is a main vertex and the other is a subdividing vertex of $X$.

If the endvertices of $P$ lie on the same edge-path of $X$, we can insert them into the appropriate edge-path of the drawings of Figure~\ref{fig:alternates}. Draw the path $P$ close to the appropriate edge-path between the endvertices of $P$ so that $P$ crosses only the edges that the edge-path cross. The resulting drawing has at most two crossings and contains a drawing of $G=G(X)\cup S$ where only edges of $G(X)$ cross.

If the endvertices of $P$ do not lie on the same $X$-edge-path, then without loss of generality $P$
connects the main vertex $b_1$ to a subdividing vertex $q$ of the $X$-edge-path  $a_1$-$b_2$.
As illustrated on 
Figure~\ref{fig:geza2},
 $X\cup P$ contains two subdivided $K_{3,3}'s$, $X$
 and $X'$, where the main vertex classes of $X'$ are
 $\{q,a_2,a_3\}$ and $\{b_1,b_2,b_3\}$.
As  $G(X)\cup S$ has exactly two cross-inducing subgraphs, it has no cross-inducing subgraph containing edges of $S$. Therefore, if $f\in S$ then $f$ must connect two vertices that are both $X$-cooperative and $X'$-cooperative. Consequently, $S$ contains no edges incident upon $a_1$, so $S\subseteq\{\gamma_{a_2a_3}$, $\gamma_{b_1b_2}$, $\gamma_{b_2b_3}$, $\gamma_{b_1b_3}\}$. 
Consider the drawing on 
Figure~\ref{fig:alternates}: insert subdividing  vertex $q$ on the $a_1$-$b_2$ edge
and connect $q$ to $b_1$ by a straight line, representing the  path $P$.
The resulting drawing has one crossing, and contains a subdrawing of $G$ in which only edges of $G(X)$ cross.

This finishes the proof of Theorem~\ref{th:geza}.

\bigskip

\noindent {\bf Remark.} AI tools were used only to check the final version of this manuscript.


\begin{thebibliography}{10}

\bibitem{bt22}
J. Bar\'at, G. T\'oth
\newblock Improvement on the Crossing Number of Crossing-Critical Graphs. 
\newblock Discrete and Computational Geometry, 67 (2022), 595--604. 
\newblock Also in: International Symposium on Graph Drawing and Network Visualization pp. 372--381, Cham: Springer International Publishing, 2020.

\bibitem{bienstock} D. Bienstock, N. Dean
\newblock Bounds for rectilinear crossing numbers
\newblock Journal of Graph Theory, Volume 17, Issue 3, July 1993, 333--348
\newblock \href{https://onlinelibrary.wiley.com/doi/10.1002/jgt.3190170308}{https://onlinelibrary.wiley.com/doi/10.1002/jgt.3190170308}

\bibitem{bdc} D. Bokal, Z. Dvo\v rák, P. Hlin\v en\'y, J. Lea\~{n}os, B. Mohar, T. Wiedera
\newblock Bounded degree conjecture holds precisely for $c$-crossing-critical graphs with $c\le 12$. 
\newblock Combinatorica, 42(2022), 701--728.
  
\bibitem{2cc} D. Bokal, B. Oporowski, R.B. Richter, G. Salazar,
\newblock Characterizing 2-crossing-critical graphs,
\newblock Advances in Applied Mathematics,  74(2016)  March, 23--208.
\newblock \href{https://doi.org/10.1016/j.aam.2015.10.003}{https://doi.org/10.1016/j.aam.2015.10.003}

\bibitem{fary} I. F\'ary, 
\newblock On straight-line representation of planar graphs
\newblock  Acta Sci. Math. (Szeged), 11: 229--233, (1948) 

\bibitem{kochol} M. Kochol, 
\newblock Construction of crossing-critical graphs, 
\newblock Discrete Math. 66 (1987), 311--313
\newblock \href{https://doi.org/10.1016/0012-365X(87)90108-7}{https://doi.org/10.1016/0012-365X(87)90108-7}

  \bibitem{kuratowski} K. Kuratowski, 
  \newblock Sur le probl\`eme des courbes gauches en topologie
  \newblock Fund. Math. (in French) 15 (1930), 271--283. 

  
  
  \bibitem{siran} J. \u{S}ir\'a\u{n}, 
  \newblock Crossing-critical edges and Kuratowski subgraphs of a graph
  \newblock JCT B 35/2 (1983), 83--92. 
\newblock \href{https://doi.org/10.1016/0095-8956(83)90064-3}{https://doi.org/10.1016/0095-8956(83)90064-3}



  \bibitem{curious} E. Czabarka, A. Helm, 
  \newblock Curious crossing-critical edges -- variations on an example of \u{S}ir\'a\u{n}
\newblock \href{https://arxiv.org/abs/2602.03717}{https://arxiv.org/abs/2602.03717}

\bibitem{wagner}
K. Wagner
\newblock Bemerkungen zum Vierfarbenproblem 
\newblock Jahresbericht der Deutschen Mathematiker-Vereinigung (1936) 46: 26–32. 1936


  
 \end{thebibliography}
 \end{document}